\documentclass[12pt,a4paper]{article}
\usepackage{amsfonts}

\usepackage[french]{babel}
\usepackage{graphicx}
\usepackage{amsmath}

\newtheorem{theorem}{Th\'{e}or\`{e}me}[section]

\newtheorem{corollary}[theorem]{Corollaire}

\newtheorem{definition}[theorem]{Definition}

\newtheorem{lemma}[theorem]{Lemme}

\newtheorem{proposition}[theorem]{Proposition}
\newtheorem{remark}[theorem]{Remarque}

\newenvironment{proof}[1][D\'{e}monstration]{\noindent\textbf{#1.} }{\hfill $\Box $}

\begin{document}

\title{Equations diff\'{e}rentielles sur les hypersurfaces de $\mathbb{P}^{4} $}
\author{Erwan Rousseau\thanks{\textit{{Adresse e-mail}: eroussea@math.uqam.ca}}}
\date{}
\maketitle

\begin{center}
{\small {D\'epartement de Math\'ematiques, Universit\'e du
Qu\'ebec \`a Montr\'eal,\\
C.P. 8888, Succursale Centre-Ville, Montr\'eal, CANADA H3C 3P8} }
\end{center}

\begin{abstract}
Dans cet article nous d\'{e}montrons que toute courbe enti\`{e}re dans une
hypersurface lisse de degr\'{e} $d\geq 97$ de $\mathbb{P}_{\mathbb{C}}^{4}$
doit satisfaire une \'{e}quation diff\'{e}rentielle alg\'{e}brique d'ordre
3. Nous donnons \'{e}galement la version logarithmique de ce
th\'{e}or\`{e}me en montrant que toute courbe enti\`{e}re dans le
compl\'{e}mentaire d'une surface lisse de degr\'{e} $d\geq 92$ de $\mathbb{P}%
_{\mathbb{C}}^{3}$ doit satisfaire une \'{e}quation diff\'{e}rentielle
alg\'{e}brique d'ordre 3. \newline

\begin{center}
\textbf{Abstract}
\end{center}

In this article we prove that every entire curve in a smooth
hypersurface of degree $d\geq 97$ in $\mathbb{P}_{\mathbb{C}}^{4}$
must satisfy an algebraic differential equation of order 3. The
logarithmic version of this result is given proving that every
entire curve in the complement of a smooth surface of degree
$d\geq 92$ in $\mathbb{P}_{\mathbb{C}}^{3}$ must satisfy an
algebraic differential equation of order 3.\newline

\noindent\textit{MSC}: 32Q45, 14F99\newline \textit{Mots-cl\'es}:
Hyperbolicit\'e au sens de Kobayashi; courbes enti\`eres;
Compl\'ementaires d'hypersurfaces projectives.
\end{abstract}

\section{Introduction}

En 1970, S.Kobayashi \cite{Ko70} a propos\'{e} la conjecture qui stipule
qu'une hypersurface g\'{e}n\'{e}rique $X$ de $\mathbb{P}_{\mathbb{C}}^{n}$
de grand degr\'{e} $d$ par rapport \`{a} $n$ est hyperbolique et que son
compl\'{e}mentaire est hyperbolique pour $d\geq 2n+1$. Dans le cas de la
dimension 2, Demailly et El Goul ont montr\'{e} \cite{DEG00}
l'hyperbolicit\'{e} des hypersurfaces tr\`{e}s g\'{e}n\'{e}riques $X\subset
\mathbb{P}_{\mathbb{C}}^{3}$ telles que $d\geq 21$ et El Goul a obtenu \cite
{E.G} que le compl\'{e}mentaire d'une courbe tr\`{e}s g\'{e}n\'{e}rique dans
$\mathbb{P}_{\mathbb{C}}^{2}$ est hyperbolique pour $d\geq 15.$ Ces
r\'{e}sultats ont \'{e}t\'{e} obtenus par l'utilisation des fibr\'{e}s de
jets de diff\'{e}rentielles de Demailly $E_{k,m}T_{X}^{\ast }$ introduits
dans \cite{De95}, et leur version logarithmique \cite{DL96}. Soit $\phi :%
\mathbb{C\rightarrow }X$ une courbe enti\`{e}re. L'id\'{e}e pour traiter les
probl\`{e}mes d'hyperbolicit\'{e} est que les sections globales des
fibr\'{e}s $E_{k,m}T_{X}^{\ast }$ s'annulant sur un diviseur ample
fournissent des \'{e}quations diff\'{e}rentielles pour la courbe $\phi .$
Dans le cas des surfaces lisses $X\subset \mathbb{P}_{\mathbb{C}}^{3},$ on
sait \cite{De95} que pour $d\geq 15$, $H^{0}(X,E_{2,m}T_{X}^{\ast }\otimes
A^{-1})\neq 0$ pour $m$ suffisamment grand$.$ Dans ce papier nous allons
montrer qu'en dimension 3, il faut au minimum chercher des op\'{e}rateurs
diff\'{e}rentiels d'ordre 3:

\begin{theorem}
\label{teo1}\textit{Soit }$X\subset \mathbb{P}^{4}$\textit{\ une
hypersurface de degr\'{e} }$d\geq 2,$ \textit{lisse. Alors: }
\begin{equation*}
H^{0}(X,E_{2,m}T_{X}^{\ast })=0.
\end{equation*}
\textit{\noindent Autrement dit, il n'y a pas de jets de diff\'{e}rentielles
d'ordre 2 d\'{e}finis globalement sur X.}
\end{theorem}

Une d\'{e}monstration simple de ce r\'{e}sultat, sans faire appel au
th\'{e}or\`{e}me de Borel-Weil-Bott \cite{DZ}, est obtenue \`{a} l'aide des
complexes de Schur.

\bigskip

R\'{e}cemment Y.T. Siu \cite{SY04} a pr\'{e}sent\'{e} une m\'{e}thode pour
produire des op\'{e}rateurs diff\'{e}rentiels en toute dimension $n$ pour un
degr\'{e} $d$ de $X$ plus grand qu'une constante d\'{e}pendant de $n.$ Si
nous nous int\'{e}ressons au degr\'{e} $d,$ le r\'{e}sultat principal de cet
article est l'obtention de sections globales des fibr\'{e}s d'op\'{e}rateurs
diff\'{e}ren\-tiels d'ordre 3 avec le r\'{e}sultat suivant

\begin{theorem}
\label{teo2}Soit $X$ une hypersurface lisse de degr\'{e} $d\geq 97$ de $%
\mathbb{P}^{4}$ et A un fibr\'{e} en droites ample, alors il y a des
sections globales de $E_{3,m}T_{X}^{\ast }\otimes A^{-1}$ pour $m$
suffisamment grand et toute courbe enti\`{e}re $f:\mathbb{C\rightarrow }X$
doit satisfaire l'\'{e}quation diff\'{e}rentielle correspondante.
\end{theorem}

\noindent et sa version logarithmique

\begin{theorem}
\label{teo3}Soit $X$ une hypersurface lisse de degr\'{e} $d\geq 92$ de $%
\mathbb{P}^{3}$ et A un fibr\'{e} en droites ample, alors il y a des
sections globales de $E_{3,m}\overline{T_{\mathbb{P}^{3}}}^{\ast }\otimes
A^{-1}$ pour $m$ suffisamment grand et toute courbe enti\`{e}re $f:\mathbb{%
C\rightarrow P}^{3}\backslash X$ doit satisfaire l'\'{e}quation
diff\'{e}rentielle correspondante.
\end{theorem}

Indiquons les points essentiels de la preuve de ces deux th\'{e}or\`{e}mes.
Elle se fonde tout d'abord sur les r\'{e}sultats que nous avons obtenus dans
\cite{Rou05} qui donnent la d\'{e}composition du gradu\'{e} du fibr\'{e} des
jets d'ordre 3 en ses repr\'{e}sentations irr\'{e}ductibles de Schur.

La cl\'{e} de la d\'{e}monstration est alors une majoration de $%
h^{2}(X,Gr^{\bullet }E_{3,m}T_{X}^{\ast }).$ Nous l'obtenons par l'\'{e}tude
de la cohomologie du fibr\'{e} vectoriel $\Gamma ^{(\lambda _{1},\lambda
_{2},\lambda _{3})}T_{X}^{\ast }$ sur $X.$ Celle-ci est faite en remontant
aux vari\'{e}t\'{e}s de drapeaux: sur celles-ci existent des fibr\'{e}s en
droites dont la cohomologie est reli\'{e}e \`{a} celle du fibr\'{e}
vectoriel consid\'{e}r\'{e} par les suites spectrales de Leray.

\section{Pr\'{e}liminaires}

\subsection{Fibr\'{e}s de jets de diff\'{e}rentielles}

Nous rappelons ici les d\'{e}finitions et propri\'{e}t\'{e}s de base des
espaces de jets construits par J.-P. Demailly dans \cite{De95}.

Soit $X$ une vari\'{e}t\'{e} complexe lisse. On part du couple $(X,V)$
o\`{u} $V\subset T_{X}$ est un sous fibr\'{e} du fibr\'{e} tangent de $X.$
On d\'{e}finit alors $X_{1}:=\mathbb{P}(V),$ et le fibr\'{e} $V_{1}\subset
T_{X_{1}}$ est d\'{e}fini par
\begin{equation*}
V_{1,(x,[v])}:=\{\xi \in T_{X_{1},(x,[v])}\text{ };\text{ }\pi _{\ast }\xi
\in \mathbb{C}v\}
\end{equation*}
o\`{u} $\pi :X_{1}\rightarrow X$ est la projection naturelle. Si l'on a un
germe de courbe $f:(\mathbb{C},0)\rightarrow (X,x)$ tangent \`{a} $V$, on
peut le relever \`{a} $X_{1}$ et l'on note le relev\'{e} $f_{[1]}.$ Dans la
suite, nous consid\'ererons le cas o\`u $V:=T_{X}$.

Par induction, on obtient une tour de vari\'{e}t\'{e}s $(X_{k},V_{k}).$ On
note $\pi _{k}:X_{k}\rightarrow X$ la projection. Consid\'{e}rons l'image
directe $\pi _{k\ast }(\mathcal{O}_{X_{k}}(m)).$ C'est un fibr\'{e}
vectoriel sur $X$ que l'on peut d\'{e}crire avec des coordonn\'{e}es
locales. Soit $z=(z_{1},...,z_{n})$ des coordonn\'{e}es centr\'{e}es en un
point $x\in X.$ Une section locale du fibr\'{e} $\pi _{k\ast }(\mathcal{O}%
_{X_{k}}(m))$ est un polyn\^{o}me
\begin{equation*}
P=\underset{\left| \alpha _{1}\right| +2\left| \alpha _{2}\right|
+...+k\left| \alpha _{k}\right| =m}{\sum }R_{\alpha }(z)dz^{\alpha
_{1}}...d^{k}z^{\alpha _{k}}
\end{equation*}
qui agit de mani\`{e}re naturelle sur les fibres du fibr\'{e} $%
J_{k}\rightarrow X$ des $k$-jets de germes de courbes dans $X$, i.e
l'ensemble des classes d'\'{e}quivalence des applications holomorphes $f:(%
\mathbb{C},0)\rightarrow (X,x)$ modulo la relation d'\'{e}quivalence
suivante: $f\sim g$ si et seulement si toutes les d\'{e}riv\'{e}es $%
f^{(j)}(0)=g^{(j)}(0)$ co\"{i}ncident pour $0\leq j\leq k$. De plus $P$ est
invariant par reparam\'{e}trisation, i.e
\begin{equation*}
P((f\circ \phi )^{\prime },...,(f\circ \phi )^{(k)})_{t}=\phi ^{\prime
}(t)^{m}P(f^{\prime },...,f^{(k)})_{\phi (t)}
\end{equation*}
pour tout $\phi $ du groupe $\mathbb{G}_{k}$ des germes de $k$%
-biholomorphismes de $(\mathbb{C},0)$. Le fibr\'{e} vectoriel $\pi _{k\ast }(%
\mathcal{O}_{X_{k}}(m))$ est not\'{e} $E_{k,m}T_{X}^{\ast }.$ En
g\'{e}n\'{e}ral, il semble difficile d'obtenir une filtration de ce
fibr\'{e}, cependant on a les cas simples suivants.

Pour $k=1,$ $E_{1,m}T_{X}^{\ast }=S^{m}T_{X}^{\ast }.$

Si $X$ est une surface, la description de $E_{2,m}T_{X}^{\ast }$ est la
suivante. On note $W=dz_{1}d^{2}z_{2}-dz_{2}d^{2}z_{1}$ le wronskien, alors
tout op\'{e}rateur diff\'{e}rentiel invariant d'ordre 2 et de degr\'{e} $m$
s'\'{e}crit
\begin{equation*}
P=\underset{\left| \alpha \right| +3k=m}{\sum }R_{\alpha ,k}(z)dz^{\alpha
}W^{k}.
\end{equation*}
Le th\'{e}or\`{e}me suivant permet de comprendre pourquoi les espaces de
jets sont utiles pour les probl\`{e}mes d'hyperbolicit\'{e}:

\bigskip

\textbf{Th\'{e}or\`{e}me \cite{GG80}, \cite{De95}. }\textit{Supposons qu'il
existe des entiers }$k,m>0$\textit{\ et un fibr\'{e} en droites ample }$L$%
\textit{\ sur }$X$\textit{\ tel que }
\begin{equation*}
H^{0}(X_{k},\mathcal{O}_{X_{k}}(m)\otimes \pi _{k}^{\ast }L^{-1})\simeq
H^{0}(X,E_{k,m}T_{X}^{\ast }\otimes L^{-1})
\end{equation*}
\textit{ait des sections non nulles }$\sigma _{1},...,\sigma _{N}.$\textit{\
Soit }$Z\subset X_{k}$\textit{\ le lieu de base de ces sections. Alors toute
courbe enti\`{e}re }$f:\mathbb{C}\rightarrow X$\textit{\ v\'{e}rifie }$%
f_{[k]}(\mathbb{C})\subset Z.$\textit{\ Autrement dit, pour tout
op\'{e}rateur diff\'{e}rentiel }$P$\textit{, }$\mathbb{G}_{k}-$\textit{%
invariant \`{a} valeurs dans }$L^{-1},$\textit{\ toute courbe enti\`{e}re }$%
f:\mathbb{C}\rightarrow X$\textit{\ v\'{e}rifie l'\'{e}quation
diff\'{e}rentielle }$P(f)=0.$

\bigskip

Rappelons qu'une vari\'{e}t\'{e} complexe, lisse et compacte est
hyperbolique s'il n'existe pas de courbes enti\`{e}res non constantes $f:%
\mathbb{C}\rightarrow X$. Ainsi, le probl\`{e}me est de produire
suffisamment d'\'{e}quations diff\'{e}rentielles alg\'{e}briquement
ind\'{e}pen\-dantes.

Dans le cas des surfaces, le r\'{e}sultat suivant donne un crit\`{e}re
num\'{e}rique pour l'existence de sections globales non nulles de $%
E_{2,m}T_{X}^{\ast }\otimes L^{-1}$

\textbf{Th\'{e}or\`{e}me \cite{De95}. }\textit{Soit X une surface
alg\'{e}brique de type g\'{e}n\'{e}ral et L un fibr\'{e} en droites ample
sur X. Alors}\textbf{\ }
\begin{equation*}
h^{0}(X,E_{2,m}T_{X}^{\ast }\otimes L^{-1})\geq \frac{m^{4}}{648}%
(13c_{1}^{2}-9c_{2})+O(m^{3}).
\end{equation*}
Pour les hypersurfaces $X\subset \mathbb{P}^{3}$ de degr\'{e} $d,$ ce
th\'{e}or\`{e}me montre l'existence de telles sections pour $d\geq 15.$

\subsection{D\'{e}composition des jets en repr\'{e}sentations
irr\'{e}duc\-tibles}

\subsubsection{Les foncteurs de Schur}

Soit $V$ est un espace vectoriel complexe de dimension finie $r$. A
l'ensemble des $r$-uplets d\'{e}croissants $(a_{1},...,a_{r})\in \mathbb{Z}%
^{r},a_{1}\geq a_{2}...\geq a_{r},$ on associe de mani\`{e}re fonctorielle
une collection d'espaces vectoriels $\Gamma ^{(a_{1},...,a_{r})}V$ qui
fournit la liste de toutes les repr\'{e}sentations polyn\^{o}miales
irr\'{e}ductibles du groupe lin\'{e}aire $Gl(V),$ \`{a} isomorphisme
pr\`{e}s. $\Gamma ^{\bullet }$ est appel\'{e} foncteur de Schur. Donnons une
description simple de ces foncteurs. Soit $\mathbb{U}_{r}=\left\{ \left(
\begin{array}{cc}
1 & \ast \\
0 & 1
\end{array}
\right) \right\} $ le groupe des matrices unipotentes triangulaires
sup\'{e}rieures $r\times r.$ Si tous les $a_{j}$ sont positifs, on
d\'{e}finit
\begin{equation*}
\Gamma ^{(a_{1},...,a_{r})}V\subset S^{a_{1}}V\otimes ...\otimes S^{a_{r}}V
\end{equation*}
comme \'{e}tant l'ensemble des polyn\^{o}mes $P(x_{1},...,x_{r})$ sur $%
(V^{\ast })^{r\text{ }}$qui sont homog\`{e}nes de degr\'{e} $a_{j}$ par
rapport \`{a} $x_{j}$ et qui sont invariants sous l'action \`{a} droite de $%
\mathbb{U}_{r}$ sur $(V^{\ast })^{r\text{ }}$ i.e tels que
\begin{equation*}
P(x_{1},...,x_{j-1},x_{j}+x_{k},x_{j+1},...,x_{r})=P(x_{1},...,x_{r})\text{
\ }\forall k<j.
\end{equation*}
Si $(a_{1},...,a_{r})$ n'est pas d\'{e}croissant alors on pose $\Gamma
^{(a_{1},...,a_{r})}V=0.$ Comme cas particuliers on retrouve les puissances
sym\'{e}triques et les puissances ext\'{e}rieures:
\begin{eqnarray*}
S^{k}V &=&\Gamma ^{(k,0,...,0)}V, \\
\wedge ^{k}V &=&\Gamma ^{(1,...,1,0,...,0)}V\text{ (avec }k\text{ indices 1),%
} \\
\det V &=&\Gamma ^{(1,...,1)}V.
\end{eqnarray*}
Les foncteurs de Schur satisfont la formule
\begin{equation*}
\Gamma ^{(a_{1}+l,...,a_{r}+l)}V=\Gamma ^{(a_{1},...,a_{r})}V\otimes (\det
V)^{l}
\end{equation*}
qui peut \^{e}tre utilis\'{e}e pour d\'{e}finir $\Gamma
^{(a_{1},...,a_{r})}V $ si l'on a des $a_{i}$ n\'{e}gatifs.

\subsubsection{Dualit\'{e} de Schur, Tableaux de Young (cf.\protect\cite{Fu.}%
)}

Faisons le lien avec les repr\'{e}sentations du groupe sym\'{e}triques. Les
repr\'{e}sen\-tations irr\'{e}ductibles du groupe sym\'{e}trique $S_{r}$
correspondent aux classes de conjugaison de $S_{r},$ i.e aux partitions $%
(l):r=l_{1}+l_{2}+...+l_{d}$ avec $l_{i}\in \mathbb{Z}$ et $l_{1}\geq
l_{2}\geq ...\geq l_{d}>0.$ La partition $(l)$ peut \^{e}tre d\'{e}crite par
un diagramme de Young avec $r$ cases et dont les longueurs des lignes sont $%
l_{1},l_{2},...,l_{d}.$ Les longueurs de ses colonnes sont $d_{j}=card\{i\in
\mathbb{Z}:l_{i}\geq j\}$ $(j=1,2,...,l_{1};l=l_{1}:$ longueur du diagramme;
$d=d_{1}:$ hauteur du diagramme; $\sum l_{i}=\sum d_{j}=r).$ Un tableau de
Young $t$ associ\'{e} \`{a} un diagramme de Young est tout simplement une
num\'{e}rotation des cases par les entiers $1,2,...,r.$ Pour un tableau de
Young fix\'{e} $t$ on introduit un idempotent $e_{t}$ de l'alg\`{e}bre $%
\mathbb{C}$.$S_{r}:$%
\begin{equation*}
e_{t}=\frac{v_{(l)}}{r!}.(\underset{q\in Q_{t}}{\sum }sgn(q).q).(\underset{%
p\in P_{t}}{\sum }p)
\end{equation*}

avec $v_{(l)}=\frac{r!}{d!}.\prod_{i=1}^{d}\frac{i!}{(l_{i}+d-i)!}.\underset{%
1\leq i<j\leq d}{\prod }(\frac{l_{i}-l_{j}}{j-i}+1)$ et les sous-groupes
\begin{eqnarray*}
P_{t} &=&\{p\in S_{r}:p\text{ pr\'{e}serve chaque ligne de }t\}, \\
Q_{t} &=&\{q\in S_{r}:q\text{ pr\'{e}serve chaque colonne de }t\}.
\end{eqnarray*}

Un tableau de Young $t$ est appel\'{e} tableau standard si sur chaque ligne
et chaque colonne les entiers sont rang\'{e}s par ordre croissant. Le nombre
de tableaux standards associ\'{e} \`{a} un diagramme de Young est \'{e}gal
\`{a} $v_{(l)}.$ Soit $D(r)$ l'ensemble de tous les tableaux standards \`{a}
$r$ cases. Alors l'identit\'{e} $1\in \mathbb{C}$.$S_{r}$ se d\'{e}compose
en
\begin{equation*}
1=\underset{t\in D(r)}{\sum }e_{t}
\end{equation*}
et ces idempotents sont orthogonaux deux \`{a} deux.

Le groupe sym\'{e}trique $S_{r}$ et donc l'alg\`{e}bre $\mathbb{C}$.$S_{r}$
agit sur $V^{\otimes r}$ par permutations des indices des \'{e}l\'{e}ments
de tenseurs:
\begin{equation*}
p(a_{1}\otimes a_{2}...\otimes a_{r})=a_{p^{-1}(1)}\otimes
a_{p^{-1}(2)}...\otimes a_{p^{-1}(r)},\text{ }\forall p\in S_{r}.
\end{equation*}
Par la d\'{e}composition pr\'{e}c\'{e}dente de l'identit\'{e} on obtient
\begin{equation*}
V^{\otimes r}=\underset{t\in D(r)}{\oplus }\Gamma ^{t}V
\end{equation*}
avec $\Gamma ^{t}V=e_{t}(V^{\otimes r}).$ Si les tableaux de Young $t,%
\widetilde{t}$ correspondent au m\^{e}me diagramme de Young, i.e \`{a} la
m\^{e}me partition $(l)$ alors $\Gamma ^{t}V$ et $\Gamma ^{\widetilde{t}}V$
sont isomorphes. D'o\`{u}
\begin{equation*}
V^{\otimes r}=\underset{(l)}{\oplus }(\Gamma ^{(l)}V)^{\oplus v_{(l)}},\text{
o\`{u} }(l)\text{ d\'{e}crit les partitions de }r.
\end{equation*}

\paragraph{Cas du fibr\'{e} cotangent}

Consid\'{e}rons une vari\'{e}t\'{e} alg\'{e}brique lisse $X$ et son
fibr\'{e} cotangent $T_{X}^{\ast }.$ Pour chaque partition $(l),$
\begin{equation*}
\dim H^{0}(X,\Gamma ^{(l)}T_{X}^{\ast })
\end{equation*}
est un invariant birationnel de la vari\'{e}t\'{e} $X$ (cf.\cite{Ma2})$.$
Dans le cas de certaines partitions particuli\`{e}res, on a par exemple le
genre cotangentiel
\begin{equation*}
\dim H^{0}(X,S^{m}T_{X}^{\ast })
\end{equation*}
(cf. \cite{Sa2}) pour $(l)=(m,0,...,0),$ ou bien le nombre de Hodge $h^{0r}$
pour $(l)=(1,..,1,0...,0)$ avec $r$ fois ''1''. Malheureusement le calcul de
ces invariants et celui des groupes de cohomologie d'ordre sup\'{e}rieur est
assez difficile (cf. chapitre 4).

\bigskip

\paragraph{Coefficient de Littlewood-Richardson}

Un diagramme de \newline
Young \textit{gauche} est le diagramme obtenu en enlevant un diagramme plus
petit d'un diagramme de Young qui le contient. Si deux diagrammes
correspondent aux partitions $\lambda =(\lambda _{1},\lambda _{2},...)$ et $%
\mu =(\mu _{1},\mu _{2},...),$ on \'{e}crit $\mu \subset \lambda $ si le
diagramme de $\mu $ est contenu dans celui de $\lambda ,$ i.e, $\mu _{i}\leq
\lambda _{i}$ pour tout $i.$ Le diagramme gauche est not\'{e} $\lambda /\mu
. $ Un tableau gauche est un diagramme gauche rempli par des entiers
positifs qui sont en croissance stricte sur les colonnes et croissance
faible sur les lignes. On d\'{e}finit \textit{le mot} d'un tableau gauche $t$
(ou mot en ligne), $w(t)$ (ou $w_{r}(t))$ en lisant les entiers de $t$ de la
gauche vers la droite et de bas en haut. Un mot $w=x_{1}x_{2}...x_{r}$ est
dit de \textit{Yamanouchi} si, quand on le lit en partant de la fin
jusqu'\`{a} n'importe quelle lettre, la suite $x_{r},x_{r-1},...,x_{s}$
contient au moins autant de ''1'' que de''2'', au moins autant de ''2'' que
de ''3'',...

Un tableau gauche $t$ est un tableau gauche de Littlewood-Richardson si son
mot $w_{r}(t)$ est de Yamanouchi. Un tableau gauche a pour contenu $\mu
=(\mu _{1},...,\mu _{l})$ si les entiers qu'il contient v\'{e}rifient : il y
a $\mu _{1}$ ''1'', $\mu _{2}$ ''2''..., $\mu _{l}$ ''l''.

\begin{definition}
(cf.\cite{Fu.})\label{d1} \textit{Le coefficient de Littlewood-Richardson} $%
c_{\lambda ,\mu }^{\nu }$ est le nombre de tableaux gauches de
Littlewood-Richardson de forme $\nu /\lambda $ de contenu $\mu .$
\end{definition}

\begin{remark}
\cite{Fu.} C'est aussi la multiplicit\'{e} de $\Gamma ^{\nu }V$ dans $\Gamma
^{\lambda }V\otimes \Gamma ^{\mu }V.$
\end{remark}

\subsubsection{Caract\`{e}res, fonctions de Schur}

On d\'{e}finit le caract\`{e}re multiplicatif de $T=%
\{(x=diag(x_{1},...,x_{r})\}\subset G$ le sous-groupe des matrices
diagonales inversibles comme \'{e}tant l'application:
\begin{eqnarray*}
\chi _{\lambda } &:&T\rightarrow \mathbb{C}^{\ast }, \\
d(t) &\rightarrow &t_{1}^{\lambda _{1}}...t_{r}^{\lambda _{r}}.
\end{eqnarray*}

\noindent

\noindent On pose: $X^{\lambda }=X_{1}^{\lambda _{1}}...X_{r}^{\lambda
_{r}}. $ On d\'{e}finit le caract\`{e}re formel $ch(E)(X)$ de $E$ comme
\'{e}tant le polyn\^{o}me:
\begin{equation*}
ch(E)(X)=\underset{\lambda }{\sum }(\dim E_{\lambda })X^{\lambda }
\end{equation*}
o\`{u} $E_{\alpha }=\{e\in E:x.e=x_{1}^{\alpha _{1}}...x_{m}^{\alpha _{m}}e$%
, $\forall $ $x\in $ $T$ $\}.$

\bigskip

\noindent Faisons le lien avec le caract\`{e}re de Chern.

\begin{proposition}
\label{p10}\textit{\noindent Soit V un fibr\'{e} vectoriel de rang r sur X. }

\noindent \textit{Notons la factorisation formelle:}
\begin{equation*}
\overset{r}{\underset{i=0}{\sum }}c_{i}(V)x^{i}=\overset{r}{\underset{i=1}{%
\Pi }}(1+a_{i}x)
\end{equation*}
\noindent $\mathit{.}$\textit{Soit }$\lambda =(\lambda _{1},...,\lambda
_{r})\in \Lambda (r,n)=\{(\lambda _{1},...,\lambda _{r})\in \mathbb{Z}%
^{r},\lambda _{1}\geq ...\geq \lambda _{r},\sum \lambda _{i}=n\}.$

\noindent \textit{Alors:}
\begin{equation*}
Ch(\Gamma ^{\lambda }V)=ch(\Gamma ^{\lambda }V)(e^{a_{1}},...,e^{a_{r}})%
\mathit{,}
\end{equation*}
\textit{o\`{u} }$Ch$\textit{\ d\'{e}signe le caract\`{e}re de Chern.}
\end{proposition}

\bigskip

\begin{proof}
Par \cite{Hi66} il suffit de le v\'{e}rifier pour $V$ somme directe de
fibr\'{e}s en droites.

\noindent Soit donc:

\noindent
\begin{equation*}
V=\overset{r}{\underset{i=1}{\oplus }}\xi _{i}.
\end{equation*}
Soit $\mathcal{U}=\{U_{i}\}_{i\in I\text{ }}$ un recouvrement ouvert pour
lequel les $\xi _{i}$ sont repr\'{e}sent\'{e}s par les cocycles $%
\{a_{ii}^{k,l}\}:$

\begin{equation*}
\noindent a_{ii}^{k,l}:U_{k}\cap U_{l}\rightarrow \mathbb{C}^{\ast }.
\end{equation*}

\noindent $V$ est donc repr\'{e}sent\'{e} par un cocycle:

\noindent $g_{k,l}(x)=\left(
\begin{array}{cccc}
a_{11}^{k,l}(x) &  &  &  \\
& ... &  &  \\
&  & ... &  \\
&  &  & a_{rr}^{k,l}(x)
\end{array}
\right) =diag(a_{11}^{k,l}(x),...,a_{rr}^{k,l}(x))$ pour $x\in U_{k}\cap
U_{l}.$

\noindent Soit :

\begin{equation*}
\noindent \rho :GL(V)\rightarrow GL(\Gamma ^{\lambda }V)
\end{equation*}
la repr\'{e}sentation associ\'{e}e \`{a} $\lambda .$

\noindent Le fibr\'{e} vectoriel $\Gamma ^{\lambda }V$ est donc
repr\'{e}sent\'{e} par le cocycle $h_{k,l}(x)=\rho (g_{k,l}(x)).$

\noindent Par le fait que $\Gamma ^{\lambda }V$ est somme de ses espaces de
poids, $\Gamma ^{\lambda }V=\underset{\mu \in \Lambda (r,n)}{\oplus }(\Gamma
^{\lambda }V)^{\mu }$ on obtient :

\noindent $h_{k,l}(x)=\rho (g_{k,l}(x))=diag(\left(
\begin{array}{ccc}
(a_{11}^{k,l})^{\mu _{1}}...(a_{rr}^{k,l})^{\mu _{r}} &  &  \\
& ... &  \\
&  & (a_{11}^{k,l})^{\mu _{1}}...(a_{rr}^{k,l})^{\mu _{r}}
\end{array}
\right) )$, matrice diagonale par blocs o\`{u} les blocs sont des matrices
diagonales de dimension $\dim (\Gamma ^{\lambda }V)^{\mu }.$

\noindent Ainsi:

\noindent
\begin{equation*}
\Gamma ^{\lambda }V=\underset{\mu \in \Lambda (r,n)}{\oplus }(\dim (\Gamma
^{\lambda }V)^{\mu })\xi _{1}^{\mu _{1}}\otimes ...\otimes \xi _{r}^{\mu
_{r}}.
\end{equation*}

\noindent Donc:

\noindent
\begin{equation*}
Ch(\Gamma ^{\lambda }V)=ch(\Gamma ^{\lambda }V)(e^{a_{1}},...,e^{a_{r}}).
\end{equation*}
\end{proof}

\bigskip

\noindent Rappelons maintenant \cite{Mart.S} qu'en caract\'{e}ristique 0, ce
qui est notre cas, on a acc\`{e}s au caract\`{e}re formel:

\noindent $ch(\Gamma ^{\lambda }V)=ch_{\lambda }=s_{\lambda }$ o\`{u} $%
s_{\lambda }$ est la fonction de Schur de type $\lambda $ d\'{e}finie par:
\begin{equation*}
s_{\lambda }=\frac{a_{\lambda +\delta }}{a_{\delta }}
\end{equation*}
\noindent o\`{u} pour $\alpha \in \Lambda (r,n)$: $a_{\alpha }(X)=\det
[(X_{i}^{\alpha _{j}})]$ et $\delta =(r-1,r-2,...,0).$

\subsubsection{D\'{e}composition des jets d'ordre 3 en dimension 3}

Dans \cite{Rou05} nous avons obtenu la d\'{e}composition du gradu\'{e} du
fibr\'{e} des jets d'ordre 3 en dimension 3:

\bigskip

\noindent \textbf{Th\'{e}or\`{e}me \cite{Rou05}. }\textit{Soit X une
vari\'{e}t\'{e} complexe de dimension 3, alors:}
\begin{equation*}
Gr^{\bullet }E_{3,m}T_{X}^{\ast }=\underset{0\leq \gamma \leq \frac{m}{5}}{%
\oplus }(\underset{\{\lambda _{1}+2\lambda _{2}+3\lambda _{3}=m-\gamma ;%
\text{ }\lambda _{i}-\lambda _{j}\geq \gamma ,\text{ }i<j\}}{\oplus }\Gamma
^{(\lambda _{1},\lambda _{2},\lambda _{3})}T_{X}^{\ast })
\end{equation*}
\noindent \textit{o\`{u}} $\Gamma $ \textit{est le foncteur de Schur.}

\bigskip

Cette d\'{e}composition a permis par Riemann-Roch un calcul de
caract\'{e}ristique d'Euler:

\bigskip

\noindent \textbf{Proposition \cite{Rou05}. }\textit{Soit X une hypersurface
lisse de degr\'{e} }$d$ de $\mathbb{P}^{4}$, alors
\begin{equation*}
\chi (X,E_{3,m}T_{X}^{\ast })=\frac{m^{9}}{81648\times 10^{6}}%
d(389d^{3}-20739d^{2}+185559d-358873)+O(m^{8})
\end{equation*}

On obtient alors la positivit\'{e} de la caract\'{e}ristique d'Euler:

\bigskip

\noindent \textbf{Corollaire \cite{Rou05}. }\textit{Pour} $d\geq 43,$ $\chi
(X,E_{3,m}T_{X}^{\ast })\sim \alpha (d)m^{9}$ \textit{avec} $\alpha (d)>0.$

\subsection{Les vari\'{e}t\'{e}s de drapeaux}

Soit $X$ une vari\'{e}t\'{e} complexe lisse de dimension 3. Notons $%
Fl(T_{X}^{\ast })$ la vari\'{e}t\'{e} des drapeaux de $T_{X}^{\ast }$ i.e
des suites de sous-espaces vectoriels embo\^{i}t\'{e}s
\begin{equation*}
D=\{0=E_{3}\subset E_{2}\subset E_{1}\subset E_{0}=T_{X,x}^{\ast }\}.
\end{equation*}

Soit $\pi :Fl(T_{X}^{\ast })\rightarrow X.$ C'est une fibration localement
triviale dont la dimension relative est: $N=1+2=3.$

Soit $\lambda =(\lambda _{1},\lambda _{2},\lambda _{3})$ une partition telle
que $\lambda _{1}>\lambda _{2}>\lambda _{3}.$ Notons $\mathcal{L}^{\lambda }$
le fibr\'{e} en droites sur $Fl(T_{X}^{\ast })$ dont la fibre au-dessus du
drapeau pr\'{e}c\'{e}dent est $\mathcal{L}_{D}^{\lambda }=\underset{i=1}{%
\overset{3}{\otimes }}\det (E_{i-1}/E_{i})^{\otimes \lambda _{i}}.$
D'apr\`{e}s le th\'{e}or\`{e}me de Bott \cite{Bot}, si $m\geq 0:$%
\begin{eqnarray*}
\pi _{\ast }(\mathcal{L}^{\lambda })^{\otimes m} &=&\Gamma ^{m\lambda
}T_{X}^{\ast }, \\
\mathcal{R}^{q}\pi _{\ast }(\mathcal{L}^{\lambda })^{\otimes m} &=&0\text{
si }q>0.
\end{eqnarray*}

Les fibr\'{e}s $\Gamma ^{m\lambda }T_{X}^{\ast }$ et $(\mathcal{L}^{\lambda
})^{\otimes m}$ ont donc m\^{e}me cohomologie.

\subsection{Le cas logarithmique}

Soit $X$ une vari\'{e}t\'{e} lisse complexe avec un diviseur \`{a}
croisements normaux $D$. On d\'{e}finit le faisceau cotangent
logarithmique
\begin{equation*}
\overline{T_{X}}^{\ast }=T_{X}^{\ast }(\log D)
\end{equation*}
comme le faisceau localement libre engendr\'{e} par $T_{X}^{\ast
}$ et les
diff\'{e}rentielles logarithmiques $\frac{ds_{j}}{s_{j}},$ o\`{u} les $%
s_{j}=0$ sont les \'{e}quations locales des composantes
irr\'{e}ductibles de $D$.

Son dual, le fibr\'{e} tangent logarithmique
\begin{equation*}
\overline{T_{X}}=T_{X}(-\log D)
\end{equation*}
est le faisceau des germes de champs de vecteurs tangents \`{a}
$D$.

En suivant \cite{DL96}, on peut d\'efinir le faisceau $O(E_{k,m}%
\overline{T_{X}}^{\ast })$ des diff\'{e}rentielles de jets
logarithmiques, i.e, le faisceau localement libre engendr\'{e} par
tous les op\'{e}rateurs polyn\^{o}miaux en les d\'{e}riv\'{e}es
d'ordre $1,2,...k$ de $f$, germe de courbe holomorphe, auxquelles
on ajoute celles de la fonction $\log (s_{j}(f))$ le long de la
j-\`{e}me composante de $D$, qui de plus sont invariants par
changement de param\'{e}trisation arbitraire.

Nous renvoyons le lecteur \`a \cite{DL96} et \cite{Rou05} pour la
g\'en\'eralisation des r\'esultats pr\'ec\'edents au cas
logarithmique.

\section{Etude de la cohomologie}

Soit $X$ une vari\'{e}t\'{e} complexe lisse de dimension 3. Nous
\'{e}tudions maintenant les fibr\'{e}s $\Gamma ^{(\lambda _{1},\lambda
_{2},\lambda _{3})}T_{X}^{\ast }$ et leur cohomologie. Pour obtenir
l'existence de suffisamment d'op\'{e}rateurs diff\'{e}rentiels, i.e de
sections globales des fibr\'{e}s $E_{k,m}T_{X}^{\ast },$ il est crucial de
contr\^{o}ler les groupes de cohomologie. En dimension 2 le point cl\'{e}
est l'utilisation d'un th\'{e}or\`{e}me d'annulation de Bogomolov \cite{Bo79}
valable sur les surfaces de type g\'{e}n\'{e}ral. En dimension 3, nous
allons voir que l'on ne peut esp\'{e}rer un tel th\'{e}or\`{e}me
d'annulation. En effet, on a peu de r\'{e}sultats explicites sur la
d\'{e}termination des groupes de cohomologie. Citons tout de m\^{e}me un
th\'{e}or\`{e}me d'annulation d\^{u} \`{a} J.-P. Demailly \cite{De95}:

\begin{theorem}
\label{t10}\textit{Soit X une vari\'{e}t\'{e} alg\'{e}brique projective de
dimension }$n$\textit{, et }$L$\textit{\ un fibr\'{e} en droites sur X.
Supposons que X est de type g\'{e}n\'{e}ral est minimal (i.e. }$K_{X}$%
\textit{\ est big et nef) et soit }$a=(a_{1},...,a_{n})\in Z^{n},a_{1}\geq
...\geq a_{n}.$\textit{\ Si l'on a L pseudoeffectif et }$\left| a\right|
=\sum a_{j}>0,$\textit{\ ou L big et }$\left| a\right| \geq 0,$\textit{\
alors}
\begin{equation*}
H^{0}(X,\Gamma ^{a}T_{X}\otimes L^{\ast })=0.
\end{equation*}
\end{theorem}

Malheureusement ce th\'{e}or\`{e}me ne nous renseigne que sur le groupe $%
H^{3}$ par la dualit\'{e} de Serre. Or, ce que nous souhaitons obtenir est
un contr\^{o}le du $H^{2}.$ Nous allons donc pr\'{e}senter dans cette
section une approche \'{e}l\'{e}mentaire pour traiter ce probl\`{e}me
utilisant des outils alg\'{e}briques (essentiellement les suites exactes),
et des outils provenant de l'analyse complexe (th\'{e}or\`{e}mes
d'annulation sous hypoth\`{e}se de positivit\'{e} des fibr\'{e}s).

\noindent Rappelons la formule:
\begin{equation*}
\Gamma ^{(\lambda _{1},\lambda _{2},\lambda _{3})}T_{X}^{\ast }=\Gamma
^{(\lambda _{1}-\lambda _{3},\lambda _{2}-\lambda _{3},0)}T_{X}^{\ast
}\otimes K_{X}^{\lambda _{3}}.
\end{equation*}

\subsection{Les jets d'ordre 1 et 2}

Montrons tout d'abord la n\'{e}cessit\'{e} d'\'{e}tudier les jets d'ordre 3.

\noindent L'absence de 1-jets d\'{e}finis globalement est bien connue:

\begin{proposition}
\textbf{(\cite{Sa2}})\label{p11}
\begin{equation*}
H^{0}(X,S^{m}T_{X}^{\ast })=0\text{\textit{\ pour }}m\geq 1.
\end{equation*}
\end{proposition}

\noindent On a \'{e}galement absence des 2-jets. En effet, on a le
r\'{e}sultat suivant:

\begin{theorem}
\label{t6}\textit{Soit }$X\subset \mathbb{P}^{4}$\textit{\ une hypersurface
de degr\'{e} }$d\geq 2,$ \textit{lisse et irr\'{e}ductible. Alors: }
\begin{equation*}
H^{0}(X,E_{2,m}T_{X}^{\ast })=0.
\end{equation*}
\textit{\noindent Autrement dit, il n'y a pas de jets de diff\'{e}rentielles
d'ordre 2 d\'{e}finis globalement sur X.}
\end{theorem}

Nous allons donner une d\'{e}monstration de ce r\'{e}sultat, sans faire
appel au th\'{e}or\`{e}me de Borel-Weil-Bott \cite{DZ}, \`{a} l'aide des
complexes de Schur qui ne sont pas utilis\'{e}s dans la preuve d'un
th\'{e}or\`{e}me plus g\'{e}n\'{e}ral de P. Br\"{u}ckmann et H.G. Rackwitz
\cite{Br} qui permet de montrer aussi ce r\'{e}sultat:

\begin{theorem}
\textit{Soit }$T$\textit{\ tableau de Young, }$d_{i}$\textit{\ le nombre de
cases de la colonne }$i$\textit{\ et }$X$\textit{\ une intersection
compl\`{e}te de dimension }$p$\textit{\ de }$\mathbb{P}^{n}.$\textit{\ Alors}
\begin{equation*}
H^{0}(X,\Gamma ^{T}T_{X}^{\ast })=0\mathit{\ }\text{\textit{si}}\mathit{\ }%
\sum_{i=1}^{n-p}d_{i}<\dim X=p.
\end{equation*}
\end{theorem}

L'outil fondamental est l'existence de suites exactes. En effet, on a la
suite exacte:
\begin{equation*}
0\rightarrow T_{X}\rightarrow T_{\mathbb{P}^{4}\left| X\right. }\rightarrow
\mathcal{O}_{X}(d)\rightarrow 0.
\end{equation*}
Donc par dualit\'{e}:
\begin{equation*}
0\rightarrow \mathcal{O}_{X}(-d)\rightarrow T_{\mathbb{P}^{4}\left| X\right.
}^{\ast }\rightarrow T_{X}^{\ast }\rightarrow 0.
\end{equation*}

Rappelons la proposition suivante qui donne une r\'{e}solution de tout
fibr\'{e} de Schur:

\begin{proposition}
(\cite{La.})\label{p8} \textit{Soit }$0\rightarrow A\rightarrow B\rightarrow
C\rightarrow 0$\textit{\ une suite exacte d'espaces vectoriels. Il y a un
complexe }$C_{\mu }^{\bullet }\rightarrow \Gamma ^{\mu }C\rightarrow 0,$%
\textit{\ dont le j-\`{e}me terme est }$C_{\mu
}^{j}=\underset{\left| \nu \right| =j,\rho }{\oplus }c_{v,\rho
}^{\mu }\Gamma ^{\nu ^{\ast }}A\otimes \Gamma ^{\rho }B$\textit{\
o\`{u} }$c_{\nu ,\rho }^{\mu }$\textit{\ est le coefficient de
Littlewood-Richardson et }$\nu ^{\ast }$ d\'{e}signe la partition
conjugu\'{e}e de $\nu .$
\end{proposition}

\noindent On applique cette proposition \`{a} la suite exacte
pr\'{e}c\'{e}dente pour obtenir une r\'{e}solution de $\Gamma
^{(b_{1},b_{2},0)}T_{X}^{\ast }:$

\bigskip

\begin{proposition}
Soit $b_{1}\geq b_{2}\geq 1.$ \textit{On a la suite exacte:}
\begin{eqnarray*}
(1) &:&0\rightarrow \mathcal{O}_{X}(-2d)\otimes \Gamma
^{(b_{1}-1,b_{2}-1,0,0)}T_{\mathbb{P}^{4}\left| X\right. }^{\ast } \\
&\rightarrow &\mathcal{O}_{X}(-d)\otimes \Gamma ^{(b_{1}-1,b_{2},0,0)}T_{%
\mathbb{P}^{4}\left| X\right. }^{\ast }\oplus \mathcal{O}_{X}(-d)\otimes
\Gamma ^{(b_{1},b_{2}-1,0,0)}T_{\mathbb{P}^{4}\left| X\right. }^{\ast } \\
&\rightarrow &\Gamma ^{(b_{1},b_{2},0,0)}T_{\mathbb{P}^{4}\left| X\right.
}^{\ast }\rightarrow \Gamma ^{(b_{1},b_{2},0)}T_{X}^{\ast }\rightarrow 0.
\end{eqnarray*}
\end{proposition}

\begin{proof}
Puisque la dimension du fibr\'{e} $\mathcal{O}_{X}(-d)$ est 1, $\nu ^{\ast }$
est de la forme $(\lambda ,0,0)$ donc les seules possibilit\'{e}s pour $\nu $
sont: $\nu =(1,0,0);\nu =(1,1,0)$ puisque $\nu \subset \mu
=(b_{1},b_{2},0,0).$

\noindent Reste la d\'{e}termination des partitions $\rho $ telles que les
coefficients de Littlewood-Richardson $c_{\nu ,\rho }^{\mu }$ soient non
nuls. Par la d\'{e}finition \ref{d1}, nous devons d\'{e}terminer les
tableaux gauches de Littlewood-Richardson de type $\mu /\nu $ de contenu $%
\rho .$ Pour $\nu =(1,0,0),$ la croissance faible sur les lignes et stricte
sur les colonnes impose que le mot $w(\mu /\nu )=(2,...2,1,..,1)$ o\`{u} il
y a $b_{2}$ ''2'' et $(b_{1}-1)$ ''1'', ou $w(\mu /\nu )=(1,2,..2,1,...1)$
o\`{u} il y a $b_{1}$ ''1'' et $(b_{2}-1)$ ''2''. Pour $\nu =(1,1,0),$ la
seule possibilit\'{e} est $w(\mu /\nu )=(2,...,2,1,...1)$ avec $(b_{2}-1)$
''2'' et $(b_{2}-1)$ ''1''. Dans tous ces cas une seule partition $\rho $
convient: respectivement $((b_{1}-1),b_{2});(b_{1},(b_{2}-1))$ et $%
(b_{1}-1,b_{2}-1).$ De plus il n'y a qu'un seul tableau gauche de
Littlewood-Richardson qui convient dans chaque cas, i.e $c_{\nu ,\rho }^{\mu
}=1.$
\end{proof}

\bigskip

\noindent On consid\`{e}re maintenant la suite exacte d'Euler:
\begin{equation*}
0\rightarrow \mathcal{O\rightarrow O}(1)^{\oplus 5}\rightarrow T_{\mathbb{P}%
^{4}}\rightarrow 0.
\end{equation*}

\noindent On lui applique la proposition \ref{p8}:

\begin{proposition}
Soit $b_{1}\geq b_{2}\geq 1.$ \textit{On a la suite exacte:}
\begin{eqnarray*}
(2) &:&0\rightarrow \mathcal{O}(b_{1}+b_{2}-2)^{\oplus
s(b_{1}-1,b_{2}-1)}\rightarrow \mathcal{O}(b_{1}+b_{2}-1)^{\oplus
(s(b_{1}-1,b_{2})+s(b_{1},b_{2}-1))} \\
&\rightarrow &\mathcal{O}(b_{1}+b_{2})^{\oplus s(b_{1},b_{2})}\rightarrow
\Gamma ^{(b_{1},b_{2},0,0)}T_{\mathbb{P}^{4}}\rightarrow 0
\end{eqnarray*}
O\`{u}:
\begin{equation*}
s(x,y)=(x-y+1)\frac{(x+2)}{2}\frac{(x+3)}{3}\frac{(x+4)}{4}(y+1)\frac{(y+2)}{%
2}\frac{(y+3)}{3}.
\end{equation*}
\end{proposition}

\begin{proof}
Le complexe de Schur nous donne la suite exacte:
\begin{eqnarray*}
0 &\rightarrow &\Gamma ^{(b_{1}-1,b_{2}-1,0,0,0)}(\mathcal{O}(1)^{\oplus
5})\rightarrow \Gamma ^{(b_{1}-1,b_{2},0,0,0)}(\mathcal{O}(1)^{\oplus
5})\oplus \Gamma ^{(b_{1},b_{2}-1,0,0,0)}(\mathcal{O}(1)^{\oplus 5}) \\
&\rightarrow &\Gamma ^{(b_{1},b_{2},0,0,0)}(\mathcal{O}(1)^{\oplus
5})\rightarrow \Gamma ^{(b_{1},b_{2},0,0)}T_{\mathbb{P}^{4}}\rightarrow 0.
\end{eqnarray*}

\noindent Il nous suffit donc de d\'{e}terminer $\Gamma ^{(a,b,0,0,0)}(%
\mathcal{O}(1)^{\oplus 5}).$ Pour \noindent $V$ somme directe de fibr\'{e}s
en droites i.e $V=\overset{r}{\underset{i=1}{\oplus }}\xi _{i}$ et $\lambda
=(\lambda _{1},...,\lambda _{r})\in \Lambda (r,n)=\{(\lambda
_{1},...,\lambda _{r})\in \mathbb{Z}^{r},\lambda _{1}\geq ...\geq \lambda
_{r},\sum \lambda _{i}=n\}$, on a
\begin{equation*}
\Gamma ^{\lambda }V=\underset{\mu \in \Lambda (r,n)}{\oplus }(\dim (\Gamma
^{\lambda }V)^{\mu })\xi _{1}^{\mu _{1}}\otimes ...\otimes \xi _{r}^{\mu
_{r}}.
\end{equation*}
Donc, ici:
\begin{equation*}
\Gamma ^{(a,b,0,0,0)}(\mathcal{O}(1)^{\oplus 5})=\mathcal{O(}a+b)^{\oplus d}
\end{equation*}
o\`{u} $d$ est le rang de $\Gamma ^{(a,b,0,0,0)}(\mathcal{O}(1)^{\oplus 5}).$
Le rang des fibr\'{e}s de Schur est connue (cf.\cite{Fu.}):

\noindent le rang de $\Gamma ^{\rho }E$ est donn\'{e} par $s_{\rho
}(1,1,...,1)$ o\`{u} $s_{\rho }$ est la fonction de Schur.

\noindent On a la propri\'{e}t\'{e} (\cite{Fu.}) :
\begin{equation*}
s_{\rho }(1,1,...,1)=\underset{i<j}{\prod }\frac{\rho _{i}-\rho _{j+j-i}}{j-i%
}.
\end{equation*}
Donc:
\begin{equation*}
d=s(a,b).
\end{equation*}
\end{proof}

\bigskip

\noindent Passons aux applications au niveau de la cohomologie en utilisant
le lemme standard:

\begin{lemma}
\label{l1}\textit{Soit E un fibr\'{e} vectoriel sur X.}

\noindent \textit{a) Supposons que l'on ait une r\'{e}solution de longueur m
}$E^{\bullet }\rightarrow E\rightarrow 0,$\textit{\ et que }$%
H^{q+j}(X,E^{j})=0$\textit{\ pour tout }$j\geq 0.$\textit{\ Alors }$%
H^{q}(X,E)=0.$

\noindent \textit{b) Supposons que l'on ait une filtration de E et notons }$%
Gr^{\bullet }E$\textit{\ la somme directe des quotients successifs de la
filtration. Si }$H^{q}(X,Gr^{\bullet }E)=0$\textit{\ pour }$q\geq 0$\textit{%
\ fix\'{e}, alors }$H^{q}(X,E)=0.$
\end{lemma}

\begin{proof}
En notant $E^{j}\overset{\phi _{j}}{\rightarrow }E^{j-1},$ on obtient les
suites exactes
\begin{eqnarray*}
0 &\rightarrow &im\text{ }\phi _{1}\rightarrow E^{0}\rightarrow E\rightarrow
0, \\
0 &\rightarrow &im\text{ }\phi _{i}\rightarrow E^{i-1}\rightarrow im\text{ }%
\phi _{i-1}\rightarrow 0,\text{ pour }2\leq i\leq m-1, \\
0 &\rightarrow &E^{m}\rightarrow E^{m-1}\rightarrow im\text{ }\phi
_{m-1}\rightarrow 0.
\end{eqnarray*}
a) d\'{e}coule imm\'{e}diatement des suites exactes longues de cohomologie.

Soit
\begin{equation*}
E\supset E^{1}\supset ...\supset E^{p}\supset ...\supset E^{m}=0
\end{equation*}
la filtration. On montre par r\'{e}currence sur $p$ que $H^{q}(X,E/E^{p})=0$
gr\^{a}ce \`{a} la suite exacte
\begin{equation*}
0\rightarrow E^{p}/E^{p+1}\rightarrow E/E^{p+1}\rightarrow
E/E^{p}\rightarrow 0.
\end{equation*}
Et b) est montr\'{e}.
\end{proof}

\bigskip

\noindent On obtient la proposition:

\begin{proposition}
\label{p12}Soit $b_{1}\geq b_{2}\geq 1.$

\noindent \textit{1) Si: }$l-b_{1}-b_{2}<0$\textit{\ alors }$H^{0}(X,\Gamma
^{(b_{1},b_{2},0,0)}T_{\mathbb{P}^{4}\left| X\right. }^{\ast }\otimes
O(l))=0.$

\noindent \textit{2) Si: }$l-b_{1}-b_{2}+1<0$\textit{\ alors }$%
H^{1}(X,\Gamma ^{(b_{1},b_{2},0,0)}T_{\mathbb{P}^{4}\left| X\right. }^{\ast
}\otimes O(l))=0.$

\noindent \textit{3) Si}$:l-b_{1}-b_{2}+2<0$\textit{\ alors }$H^{2}(X,\Gamma
^{(b_{1},b_{2},0,0)}T_{\mathbb{P}^{4}\left| X\right. }^{\ast }\otimes
O(l))=0.$

\noindent \textit{4) Si : }$b_{1}+b_{2}-l+(d-5)<0$\textit{\ alors }$%
H^{3}(X,\Gamma ^{(b_{1},b_{2},0,0)}T_{\mathbb{P}^{4}\left| X\right. }^{\ast
}\otimes O(l))=0.$
\end{proposition}

\begin{proof}
On applique la partie a) du lemme \ref{l1} pr\'{e}c\'{e}dent \`{a} la suite
exacte (2):
\begin{eqnarray*}
H^{0}(\mathbb{P}^{4},\Gamma ^{(b_{1},b_{2},0,0)}T_{\mathbb{P}^{4}}\otimes
O(p)) &=&0\text{ pour }p+b_{1}+b_{2}<0, \\
H^{1}(\mathbb{P}^{4},\Gamma ^{(b_{1},b_{2},0,0)}T_{\mathbb{P}^{4}}\otimes
O(p)) &=&0, \\
H^{2}(\mathbb{P}^{4},\Gamma ^{(b_{1},b_{2},0,0)}T_{\mathbb{P}^{4}}\otimes
O(p)) &=&0\text{ pour }-3-b_{1}-b_{2}-p<0, \\
H^{3}(\mathbb{P}^{4},\Gamma ^{(b_{1},b_{2},0,0)}T_{\mathbb{P}^{4}}\otimes
O(p)) &=&0\text{ pour }-4-b_{1}-b_{2}-p<0, \\
H^{4}(\mathbb{P}^{4},\Gamma ^{(b_{1},b_{2},0,0)}T_{\mathbb{P}^{4}}\otimes
O(p)) &=&0\text{ pour }-5-b_{1}-b_{2}-p<0.
\end{eqnarray*}
Par la dualit\'{e} de Serre
\begin{equation*}
H^{q}(X,\Gamma ^{(b_{1},b_{2},0,0)}T_{\mathbb{P}^{4}\left| X\right. }^{\ast
}\otimes \mathcal{O}(l))=H^{3-q}(X,\Gamma ^{(b_{1},b_{2},0,0)}T_{\mathbb{P}%
^{4}\left| X\right. }\otimes \mathcal{O}(d-5-l)).
\end{equation*}

\noindent Montrer 1) revient \`{a} voir pour quelles conditions
\begin{equation*}
H^{3}(X,\Gamma ^{(b_{1},b_{2},0,0)}T_{\mathbb{P}^{4}/X}\otimes \mathcal{O}%
(d-5-l))=0.
\end{equation*}
On a la suite exacte:
\begin{eqnarray*}
0 &\rightarrow &\Gamma ^{(b_{1},b_{2},0,0)}T_{\mathbb{P}^{4}}\otimes
\mathcal{O}(-5-l)\rightarrow \Gamma ^{(b_{1},b_{2},0,0)}T_{\mathbb{P}%
^{4}}\otimes \mathcal{O}(d-5-l) \\
&\rightarrow &\Gamma ^{(b_{1},b_{2},0,0)}T_{\mathbb{P}^{4}\left| X\right.
}\otimes \mathcal{O}(d-5-l)\rightarrow 0.
\end{eqnarray*}

\noindent Donc: $H^{3}(X,\Gamma ^{(b_{1},b_{2},0,0)}T_{\mathbb{P}^{4}\left|
X\right. }\otimes \mathcal{O}(d-5-l))=0$ pour
\begin{eqnarray*}
H^{3}(\mathbb{P}^{4},\Gamma ^{(b_{1},b_{2},0,0)}T_{\mathbb{P}^{4}}\otimes
O(d-5-l)) &=&0, \\
H^{4}(\mathbb{P}^{4},\Gamma ^{(b_{1},b_{2},0,0)}T_{\mathbb{P}^{4}}\otimes
O(-5-l)) &=&0.
\end{eqnarray*}
Ce qui est vrai par ce qui pr\'{e}c\`{e}de pour
\begin{eqnarray*}
-4-b_{1}-b_{2}-(d-5-l) &<&0, \\
-5-b_{1}-b_{2}-(-5-l) &<&0.
\end{eqnarray*}
Et 1) est montr\'{e}.

\noindent Montrer 2) revient \`{a} voir pour quelles conditions
\begin{equation*}
H^{2}(X,\Gamma ^{(b_{1},b_{2},0,0)}T_{\mathbb{P}^{4}/X}\otimes \mathcal{O}%
(d-5-l))=0.
\end{equation*}
Ceci est le cas pour
\begin{eqnarray*}
H^{2}(\mathbb{P}^{4},\Gamma ^{(b_{1},b_{2},0,0)}T_{\mathbb{P}^{4}}\otimes
O(d-5-l)) &=&0, \\
H^{3}(\mathbb{P}^{4},\Gamma ^{(b_{1},b_{2},0,0)}T_{\mathbb{P}^{4}}\otimes
O(-5-l)) &=&0.
\end{eqnarray*}
Ce qui est v\'{e}rifi\'{e} pour
\begin{eqnarray*}
-3-b_{1}-b_{2}-(d-5-l) &<&0, \\
-4-b_{1}-b_{2}-(-5-l) &<&0.
\end{eqnarray*}
Et 2) est montr\'{e}.

\noindent Montrer 3) revient \`{a} voir pour quelles conditions
\begin{equation*}
H^{1}(X,\Gamma ^{(b_{1},b_{2},0,0)}T_{\mathbb{P}^{4}\left| X\right. }\otimes
\mathcal{O}(d-5-l))=0.
\end{equation*}
Ceci est le cas pour
\begin{eqnarray*}
H^{1}(\mathbb{P}^{4},\Gamma ^{(b_{1},b_{2},0,0)}T_{\mathbb{P}^{4}}\otimes
O(d-5-l)) &=&0, \\
H^{2}(\mathbb{P}^{4},\Gamma ^{(b_{1},b_{2},0,0)}T_{\mathbb{P}^{4}}\otimes
O(-5-l)) &=&0.
\end{eqnarray*}
Ce qui est v\'{e}rifi\'{e} pour
\begin{equation*}
-3-b_{1}-b_{2}-(-5-l)<0.
\end{equation*}
Et 3) est montr\'{e}.

\noindent Montrer 4) revient \`{a} voir pour quelles conditions
\begin{equation*}
H^{0}(X,\Gamma ^{(b_{1},b_{2},0,0)}T_{\mathbb{P}^{4}\left| X\right. }\otimes
\mathcal{O}(d-5-l))=0.
\end{equation*}
Ceci est le cas pour
\begin{eqnarray*}
H^{0}(\mathbb{P}^{4},\Gamma ^{(b_{1},b_{2},0,0)}T_{\mathbb{P}^{4}}\otimes
O(d-5-l)) &=&0, \\
H^{1}(\mathbb{P}^{4},\Gamma ^{(b_{1},b_{2},0,0)}T_{\mathbb{P}^{4}}\otimes
O(-5-l)) &=&0.
\end{eqnarray*}
Ce qui est v\'{e}rifi\'{e} pour
\begin{equation*}
(d-5-l)+b_{1}+b_{2}<0.
\end{equation*}
Et 4) est montr\'{e}.
\end{proof}

\bigskip

\noindent On en d\'{e}duit donc:

\begin{proposition}
\label{p4}Soit $b_{1}\geq b_{2}\geq 1$ , $d\geq 2.$%
\begin{equation*}
H^{0}(X,\Gamma ^{(b_{1},b_{2},0)}T_{X}^{\ast }\otimes \mathcal{O}(l))=0%
\mathit{\ }\text{\textit{pour }}l-b_{1}-b_{2}<0.
\end{equation*}
\end{proposition}

\begin{proof}
On applique la partie a) du lemme \ref{l1} \`{a} la suite exacte (1). $%
H^{0}(X,\Gamma ^{(b_{1},b_{2},0)}T_{X}^{\ast }\otimes \mathcal{O}(l))=0$
pour
\begin{eqnarray*}
H^{0}(X,\Gamma ^{(b_{1},b_{2},0,0)}T_{\mathbb{P}^{4}\left| X\right. }^{\ast
}\otimes O(l)) &=&0, \\
H^{1}(X,\mathcal{O}_{X}(-d)\otimes \Gamma ^{(b_{1}-1,b_{2},0,0)}T_{\mathbb{P}%
^{4}/X}^{\ast }\otimes O(l))&=&0, \\
H^{1}(X,\mathcal{O}_{X}(-d)\otimes \Gamma ^{(b_{1},b_{2}-1,0,0)}T_{\mathbb{P}%
^{4}\left| X\right. }^{\ast }\otimes O(l)) &=&0, \\
H^{2}(X,\mathcal{O}_{X}(-2d)\otimes \Gamma ^{(b_{1}-1,b_{2}-1,0,0)}T_{%
\mathbb{P}^{4}\left| X\right. }^{\ast }\otimes O(l)) &=&0.
\end{eqnarray*}
Ce qui est v\'{e}rifi\'{e} par la proposition \ref{p12} pour
\begin{eqnarray*}
l-b_{1}-b_{2} &<&0, \\
l-d-b_{1}-b_{2}+2 &<&0, \\
l-2d-b_{1}-b_{2}+4 &<&0.
\end{eqnarray*}
Et la proposition est montr\'{e}e.
\end{proof}

\bigskip

\noindent On peut maintenant d\'{e}montrer le th\'{e}or\`{e}me \ref{teo1}
annonc\'{e}:

\bigskip

\begin{proof}
On sait que : $Gr^{\bullet }E_{2,m}T_{X}^{\ast }=\underset{\lambda
_{1}+2\lambda _{2}=m}{\oplus }\Gamma ^{(\lambda _{1},\lambda
_{2},0)}T_{X}^{\ast }.$ On applique la proposition pr\'{e}c\'{e}dente \ref
{p4} et la proposition \ref{p11} qui nous donnent $H^{0}(X,Gr^{\bullet
}E_{2,m}T_{X}^{\ast })=0.$ Donc, par la partie b) du lemme \ref{l1},
\begin{equation*}
H^{0}(X,E_{2,m}T_{X}^{\ast })=0.
\end{equation*}
\end{proof}

\bigskip

\noindent La prochaine section va montrer que nous pouvons aussi utiliser
les th\'{e}or\`{e}mes d'annulation classiques.

\subsection{Premiers r\'{e}sultats sur les jets d'ordre 3}

Faisons tout d'abord quelques rappels. Soit $E$ un fibr\'{e} vectoriel
hermitien de rang $r$ sur une vari\'{e}t\'{e} complexe compacte $X$ de
dimension $n.$ On note par $C_{p,q}^{\infty }(E)$ l'espace des formes
diff\'{e}rentielles $C^{\infty }$ de type $(p,q)$ sur $X$ \`{a} valeurs dans
$E$ et par
\begin{equation*}
D_{E}=D_{E}^{\prime }+D_{E}^{\prime \prime }:C_{p,q}^{\infty }(E)\rightarrow
C_{p+1,q}^{\infty }(E)\oplus C_{p,q+1}^{\infty }(E)
\end{equation*}
la connection de Chern de $E$. Soient $(x_{1},...,x_{n})$ des
coordonn\'{e}es holomorphes sur $X$ et $(e_{1},...,e_{n})$ un rep\`{e}re
orthonormal mobile $C^{\infty }$ de $E.$ Le tenseur de courbure de Chern $%
c(E)$ est d\'{e}fini par $D_{E}^{2}=c(E)\wedge \bullet $ et peut
s'\'{e}crire
\begin{equation*}
c(E)=\underset{i,j,\lambda ,\mu }{\sum }c_{ij\lambda \mu }dx_{i}d\overline{%
x_{j}}\otimes e_{\lambda }^{\ast }\otimes e_{\mu },1\leq i,j\leq n,1\leq
\lambda ,\mu \leq r.
\end{equation*}
Le tenseur de courbure $ic(E)$ est en fait une $(1,1)-$forme \`{a} valeur
dans le fibr\'{e} $Herm(E,E)$ des endomorphismes hermitiens de $E$, i.e $%
c_{ij\lambda \mu }=\overline{c}_{ji\mu \lambda };$ ainsi $ic(E)$ peut
\^{e}tre identifi\'{e} avec une forme hermitienne sur $T_{X}\otimes E.$

Rappelons qu'un fibr\'{e} vectoriel $E$ est positif (respectivement
semi-positif) au sens de Griffiths si on peut munir $E$ d'une m\'{e}trique
hermitienne telle qu'en tout point $x\in X:$

$ic(E)_{x}(\zeta \otimes v,\zeta \otimes v)=\underset{i,j,\lambda ,\mu }{%
\sum }c_{ij\lambda \mu }(x)\zeta _{i}\overline{\zeta _{j}}v_{\lambda }%
\overline{v_{\mu }}$ $>0,$ resp. $\geq 0;$ o\`{u} $ic(E)$ est le tenseur de
courbure, $\zeta =\sum \zeta _{i}\frac{\partial }{\partial z_{i}}\in
T_{X},v=\sum v_{\lambda }e_{\lambda }\in E_{x}.$ Rappelons \'{e}galement que
tout fibr\'{e} engendr\'{e} par ses sections, i.e tel que l'application $%
H^{0}(X,E)\rightarrow E_{x}$ est surjective, est semi-positif (cf.\cite{De87}%
)

Nous allons utiliser un th\'{e}or\`{e}me d'annulation d\^{u} \`{a} J.-P.
Demailly \cite{De87}:

\begin{theorem}
\textit{Soit X une vari\'{e}t\'{e} complexe de dimension n et L un fibr\'{e}
en droites sur X; E un fibr\'{e} vectoriel de rang r. Supposons }$E>0$%
\textit{\ et }$L\geq 0$\textit{, ou }$E\geq 0$\textit{\ et }$L>0.$\textit{\
Soit }$h\in \{1,...,r-1\}$\textit{\ et }$\Gamma ^{a}E$\textit{\ le fibr\'{e}
de Schur de poids }$a\in Z^{r},$\textit{\ avec }$a_{1}\geq a_{2}\geq ...\geq
a_{h}>a_{h+1}=...=a_{r}=0.$

\noindent\textit{Alors pour }$q\geq 1,$\textit{\ }$H^{n,q}(X,\Gamma
^{a}E\otimes (\det E)^{l}\otimes L)=0$\textit{\ pour }$l\geq h.$
\end{theorem}

Nous allons montrer le r\'{e}sultat suivant:

\begin{theorem}
\textit{Soit }$X\subset \mathbb{P}^{4}$\textit{\ une hypersurface lisse et
irr\'{e}ductible de degr\'{e} d.}

\noindent\textit{Alors pour }$q\geq 1,$\textit{\ }
\begin{equation*}
H^{q}(X,\Gamma ^{(a_{1},a_{2},a_{3})}T_{X}^{\ast })=0\mathit{\ }\text{%
\textit{pour}}\mathit{\ }a_{3}(d-1)>2(a_{1}+a_{2})+3(d-1).
\end{equation*}
\end{theorem}

\begin{proof}
Dans notre situation on a $X\subset \mathbb{P}^{4}$ une hypersurface lisse
et irr\'{e}ductible de degr\'{e} d. On a $T_{\mathbb{P}^{4}}^{\ast }\otimes
\mathcal{O}(2)$ qui est engendr\'{e} par ses sections$,$ donc $T_{X}^{\ast
}\otimes \mathcal{O}_{X}(2),$ quotient de $T_{\mathbb{P}^{4}}^{\ast }\otimes
\mathcal{O}(2),$ est semi-positif. Appliquons le th\'{e}or\`{e}me
pr\'{e}c\'{e}dent \`{a} $E=T_{X}^{\ast }\otimes \mathcal{O}_{X}(2)\geq 0:$%
\begin{equation*}
\Gamma ^{(a_{1}-a_{3},a_{2}-a_{3},0)}E=\Gamma
^{(a_{1}-a_{3},a_{2}-a_{3},0)}T_{X}^{\ast }\otimes \mathcal{O}%
_{X}(2(a_{1}-a_{3}+a_{2}-a_{3}))
\end{equation*}
et
\begin{equation*}
\det E=K_{X}\otimes \mathcal{O}_{X}(6).
\end{equation*}
\begin{equation*}
\Gamma ^{(a_{1}-a_{3},a_{2}-a_{3},0)}E\otimes (\det E)^{2}=\Gamma
^{(a_{1}-a_{3},a_{2}-a_{3},0)}T_{X}^{\ast }\otimes \mathcal{O}%
_{X}(2(a_{1}+a_{2})-4a_{3}+2(d+1)).
\end{equation*}
Donc:
\begin{equation*}
H^{q}(X,\Gamma ^{(a_{1},a_{2},a_{3})}T_{X}^{\ast })=H^{3,q}(X,\Gamma
^{(a_{1}-a_{3},a_{2}-a_{3},0)}T_{X}^{\ast }\otimes K_{X}^{(a_{3}-1)})=0
\end{equation*}
pour
\begin{equation*}
(a_{3}-1)(d-5)>2(a_{1}+a_{2})-4a_{3}+2(d+1)
\end{equation*}
i.e
\begin{equation*}
a_{3}(d-1)>2(a_{1}+a_{2})+3(d-1).
\end{equation*}
\end{proof}

\begin{remark}
1) l'\'{e}tude des suites exactes d\'{e}crites pr\'{e}c\'{e}demment par les
complexes de Schur fournissent des r\'{e}sultats \'{e}quivalents en
utilisant le th\'{e}or\`{e}me de Borel-Weil-Bott \cite{DZ} comme l'ont
montr\'{e} \cite{Ma2} et \cite{Br}. Ainsi, l'utilisation du Th\'{e}or\`{e}me
7 de \cite{Br} \ fournit le r\'{e}sultat:
\begin{equation*}
\text{pour }q\geq 1,\text{ }H^{q}(X,\Gamma ^{(a_{1},a_{2},a_{3})}T_{X}^{\ast
})=0\mathit{\ }\text{\textit{pour}}\mathit{\ }a_{3}(d-1)>2(a_{1}+a_{2})+3d-8.
\end{equation*}
2) On obtient donc qu'\`{a} $m$ fix\'{e}, pour $d$ suffisamment grand on a
''peu'' en proportion de partitions qui donnent une contribution non nulle
du $H^{2}.$
\end{remark}

\bigskip

\noindent On obtient donc le corollaire suivant:

\begin{corollary}
\textit{Pour} $a_{3}(d-1)>2(a_{1}+a_{2})+3(d-1),$
\begin{equation*}
h^{0}(X,\Gamma ^{(a_{1},a_{2},a_{3})}T_{X}^{\ast })=\chi (X,\Gamma
^{(a_{1},a_{2},a_{3})}T_{X}^{\ast }).
\end{equation*}
\end{corollary}

\begin{remark}
Contrairement au cas des jets d'ordre 2 en dimension 2, on ne peut
esp\'{e}rer avoir $H^{2}(X,Gr^{\bullet }E_{3,m}T_{X}^{\ast })=0$ car pour
tout m suffisamment grand, il existe $H^{2}(X,\Gamma ^{(\lambda _{1},\lambda
_{2},\lambda _{3})}T_{X}^{\ast })\neq 0.$
\end{remark}

\noindent En effet:

\begin{proposition}
\textit{Soit }$X\subset \mathbb{P}^{4}$\textit{\ une hypersurface lisse et
irr\'{e}ductible de degr\'{e} }$d\geq 6.$

\noindent\textit{Alors}
\begin{equation*}
h^{2}(X,S^{m}T_{X}^{\ast })\underset{+\infty }{\sim }(-\frac{7}{24}d+\frac{1%
}{8}d^{2})m^{5}.
\end{equation*}
\end{proposition}

\begin{proof}
On sait par la proposition \ref{p11} que
\begin{equation*}
h^{0}(X,S^{m}T_{X}^{\ast })=0.
\end{equation*}
De plus $H^{3}(X,S^{m}T_{X}^{\ast })=H^{0}(X,S^{m}T_{X}\otimes K_{X})=0$
pour $m>6$ par le th\'{e}or\`{e}me \ref{t10} car $X$ est de type
g\'{e}n\'{e}ral $(d\geq 6)$ et
\begin{equation*}
S^{m}T_{X}\otimes K_{X}=\Gamma ^{(m-2,-2,-2)}T_{X}\otimes K_{X}^{-1}.
\end{equation*}
On a les suites exactes:
\begin{equation*}
(1):0\rightarrow S^{m}T_{\mathbb{P}^{4}\left| X\right. }^{\ast }\rightarrow
\underset{\binom{4+m}{m}}{\oplus }\mathcal{O}(-m)\rightarrow \underset{%
\binom{4+m-1}{m-1}}{\oplus }\mathcal{O}(1-m)\rightarrow 0
\end{equation*}
\begin{equation*}
(2):0\rightarrow S^{m-1}T_{\mathbb{P}^{4}\left| X\right. }^{\ast }\otimes
\mathcal{O}(-d)\rightarrow S^{m}T_{\mathbb{P}^{4}\left| X\right. }^{\ast
}\rightarrow S^{m}T_{X}^{\ast }\rightarrow 0.
\end{equation*}
De (1) il vient $H^{2}(X,S^{m}T_{\mathbb{P}^{4}\left| X\right. }^{\ast
}\otimes \mathcal{O}(l))=0$ pour $m>0$ et $l\in \mathbb{Z}$ et $%
H^{1}(X,S^{m}T_{\mathbb{P}^{4}\left| X\right. }^{\ast })=0$ pour $m\geq 2.$
Donc par (2) $H^{1}(X,S^{m}T_{X}^{\ast })=0$ pour $m\geq 2.$ Finalement
\begin{equation*}
\chi (X,S^{m}T_{X}^{\ast })=h^{2}(X,S^{m}T_{X}^{\ast })\text{ pour }m>6.
\end{equation*}
On conclut gr\^{a}ce au calcul explicite par Riemann-Roch.
\end{proof}

\section{Les jets d'ordre 3 en dimension 3}

Nous allons montrer le th\'{e}or\`{e}me

\begin{theorem}
\label{t2}Soit $X$ une hypersurface lisse de degr\'{e} $d$ de $\mathbb{P}%
^{4} $, alors
\begin{equation*}
h^{2}(X,Gr^{\bullet }E_{3,m}T_{X}^{\ast })\leq Cd(d+13)m^{9}+O(m^{8})
\end{equation*}
o\`{u} C est une constante.
\end{theorem}

La preuve s'inspire de la d\'{e}monstration alg\'{e}brique \cite{Ang} des
in\'{e}galit\'{e}s de Morse de Demailly \cite{De96} qui stipulent:

\begin{theorem}
\label{t1}Soit $L=F-G$ un fibr\'{e} en droites sur une vari\'{e}t\'{e}
compacte K\"{a}hler X o\`{u} F et G sont des fibr\'{e}s en droites nef.
Alors pour $0\leq q\leq n=\dim X$%
\begin{equation*}
h^{q}(X,L^{\otimes k})\leq \frac{k^{n}}{(n-q)!q!}F^{n-q}.G^{q}+o(k^{n}).
\end{equation*}
\end{theorem}

Montrons tout d'abord la proposition

\begin{proposition}
\label{p2}Soit $\lambda =(\lambda _{1},\lambda _{2},\lambda _{3})$ une
partition telle que $\lambda _{1}>\lambda _{2}>\lambda _{3}$ et $\left|
\lambda \right| =\sum \lambda _{i}>4(d-5)+18.$ Alors :
\begin{equation*}
h^{2}(Fl(T_{X}^{\ast }),\mathcal{L}^{\lambda })=h^{2}(X,\Gamma ^{\lambda
}T_{X}^{\ast })\leq g(\lambda )d(d+13)+q(\lambda )
\end{equation*}
o\`{u} $g(\lambda )=\frac{3\left| \lambda \right| ^{3}}{2}\underset{\lambda
_{i}>\lambda _{j}}{\prod }(\lambda _{i}-\lambda _{j})$ et de plus $q$ est un
polyn\^{o}me en $\lambda $ de composantes homog\`{e}nes de plus haut
degr\'{e} 5.
\end{proposition}

\begin{proof}
On a
\begin{equation*}
\mathcal{L}^{\lambda }=(\mathcal{L}^{\lambda }\otimes \pi ^{\ast }\mathcal{O}%
_{X}(3\left| \lambda \right| ))\otimes (\pi ^{\ast }\mathcal{O}_{X}(3\left|
\lambda \right| ))^{-1}=F\otimes G^{-1},
\end{equation*}
avec $F=\mathcal{L}^{\lambda }\otimes \pi ^{\ast }\mathcal{O}_{X}(3\left|
\lambda \right| )$, $G=\pi ^{\ast }\mathcal{O}_{X}(3\left| \lambda \right|
). $ $\mathcal{L}^{\lambda }\otimes \pi ^{\ast }\mathcal{O}_{X}(3\left|
\lambda \right| )$ est positif. En effet, on a la propri\'{e}t\'{e}
g\'{e}n\'{e}rale \cite{De87} que si $E$ est un fibr\'{e} vectoriel
semi-positif i.e $E\geq 0$ alors le fibr\'{e} en droites correspondant $%
\mathcal{L(}E)^{\lambda }$ est aussi semi-positif. Ici, $E=T_{X}^{\ast
}\otimes \mathcal{O}_{X}(2)$ est semi-positif et
\begin{equation*}
\mathcal{L(}E)^{\lambda }\simeq \mathcal{L}^{\lambda }\otimes \pi ^{\ast }%
\mathcal{O}_{X}(2\left| \lambda \right| )\geq 0,
\end{equation*}
donc
\begin{equation*}
\mathcal{L}^{\lambda }\otimes \pi ^{\ast }\mathcal{O}_{X}(3\left| \lambda
\right| )>0.
\end{equation*}

Soit $Y=Fl(T_{X}^{\ast }).$ Tout d'abord montrons que $H^{i}(Y,F)=0$ pour
tout $i\geq 1$ et $\lambda $ telle que $\left| \lambda \right| =\sum \lambda
_{i}>4(d-5)+18.$ Pour cela nous utilisons le th\'{e}or\`{e}me d'annulation
de Kodaira qui stipule que pour tout fibr\'{e} en droites $A$ ample sur une
vari\'{e}t\'{e} projective $Z$ complexe $H^{i}(Z,K_{Z}\otimes A)=0$ pour $%
i>0.$ En effet, regardons \`{a} quelles conditions
\begin{equation*}
F\otimes K_{Y}^{-1}>0.
\end{equation*}
Rappelons \cite{Ma} que
\begin{equation*}
K_{Y}=\mathcal{L}^{-(5,3,1)}\otimes \pi ^{\ast }(K_{X}\otimes \det
(T_{X}^{\ast })^{\otimes 3})=\mathcal{L}^{-(5,3,1)}\otimes \pi ^{\ast }%
\mathcal{O}_{X}(4(d-5)).
\end{equation*}
Donc
\begin{equation*}
F\otimes K_{Y}^{-1}=\mathcal{L}^{\lambda +(5,3,1)}\otimes \pi ^{\ast }%
\mathcal{O}_{X}(3\left| \lambda \right| -4(d-5)).
\end{equation*}
Or on a
\begin{equation*}
\mathcal{L}^{\lambda +(5,3,1)}\otimes \pi ^{\ast }\mathcal{O}_{X}(2\left|
\lambda +(5,3,1)\right| )=\mathcal{L}^{\lambda +(5,3,1)}\otimes \pi ^{\ast }%
\mathcal{O}_{X}(2\left| \lambda \right| +18)\geq 0.
\end{equation*}
Par cons\'{e}quent $F\otimes K_{Y}^{-1}>0$ si
\begin{equation*}
3\left| \lambda \right| -4(d-5)>2\left| \lambda \right| +18
\end{equation*}
c'est-\`{a}-dire
\begin{equation*}
\left| \lambda \right| >4(d-5)+18.
\end{equation*}
Prenons un diviseur $D=\pi ^{\ast }E_{1}\in \left| G\right| ,$ lisse et
irr\'{e}ductible. On a la suite exacte:
\begin{equation*}
0\rightarrow \mathcal{O}_{Y}(F\otimes G^{-1})\rightarrow \mathcal{O}%
_{Y}(F)\rightarrow \mathcal{O}_{D}(F)\rightarrow 0.
\end{equation*}
donc la suite exacte longue en cohomologie:
\begin{equation*}
0=H^{1}(Y,\mathcal{O}_{Y}(F))\rightarrow H^{1}(D,\mathcal{O}%
_{D}(F))\rightarrow H^{2}(Y,\mathcal{O}_{Y}(F\otimes G^{-1})\rightarrow
H^{2}(Y,\mathcal{O}_{Y}(F))=0.
\end{equation*}
Donc
\begin{equation*}
h^{2}(Y,\mathcal{O}_{Y}(F\otimes G^{-1}))=h^{1}(D,\mathcal{O}_{D}(F)).
\end{equation*}
Prenons un deuxi\`{e}me diviseur $D^{\prime }=\pi ^{\ast }E_{2}\in \left|
G\right| ,$ lisse et irr\'{e}ductible, rencontrant $D$ proprement. Soit $%
Z=D\cap D^{\prime }$, $F^{\prime }=F\otimes G$ et $E_{3}=E_{1}\cap E_{2}.$
On a la suite exacte
\begin{equation*}
0\rightarrow \mathcal{O}_{D}(F^{\prime }\otimes G^{-1})\rightarrow \mathcal{O%
}_{D}(F^{\prime })\rightarrow \mathcal{O}_{Z}(F^{\prime })\rightarrow 0.
\end{equation*}
Par adjonction
\begin{equation*}
K_{D}=(K_{Y})_{|D}\otimes \mathcal{O}_{D}(D)
\end{equation*}
donc
\begin{equation*}
F_{|D}^{\prime }\otimes K_{D}^{-1}=(F\otimes K_{Y}^{-1})_{|D}>0.
\end{equation*}
Ainsi
\begin{equation*}
h^{1}(D,\mathcal{O}_{D}(F))\leq h^{0}(Z,\mathcal{O}_{Z}(F^{\prime
}))=h^{0}(Z,\mathcal{O}_{Z}(F\otimes G))\leq h^{0}(Z,\mathcal{O}%
_{Z}(F\otimes G^{2})).
\end{equation*}
Or comme pr\'{e}c\'{e}demment
\begin{equation*}
\mathcal{O}_{Z}(F\otimes G^{2})\otimes K_{Z}^{-1}=(F\otimes
K_{Y}^{-1})_{|Z}>0
\end{equation*}
donc
\begin{equation*}
h^{0}(Z,\mathcal{O}_{Z}(F\otimes G^{2}))=\chi (Z,\mathcal{O}_{Z}(F\otimes
G^{2})).
\end{equation*}
On a
\begin{equation*}
\chi (Z,\mathcal{O}_{Z}(F\otimes G^{2}))=\chi (E_{3},\Gamma ^{\lambda
}T_{X|E_{3}}^{\ast }\otimes \mathcal{O}_{E_{3}}(9\left| \lambda \right| )).
\end{equation*}
Par Riemann-Roch et la proposition \ref{p10}, on sait explicitement calculer
\begin{equation*}
\chi (X,\Gamma ^{\lambda }T_{X}^{\ast }\otimes \mathcal{O}_{X}(t)).
\end{equation*}
On a les suites exactes
\begin{eqnarray*}
(1)\text{ }0 &\rightarrow &\Gamma ^{\lambda }T_{X}^{\ast }\otimes \mathcal{O}%
_{X}(t-E_{1})\rightarrow \Gamma ^{\lambda }T_{X}^{\ast }\otimes \mathcal{O}%
_{X}(t)\rightarrow \Gamma ^{\lambda }T_{X|E_{1}}^{\ast }\otimes \mathcal{O}%
_{E_{1}}(t)\rightarrow 0 \\
(2)\text{ }0 &\rightarrow &\Gamma ^{\lambda }T_{X|E_{1}}^{\ast }\otimes
\mathcal{O}_{E}(t-E_{3})\rightarrow \Gamma ^{\lambda }T_{X|E_{1}}^{\ast
}\otimes \mathcal{O}_{E_{1}}(t)\rightarrow \Gamma ^{\lambda
}T_{X|E_{3}}^{\ast }\otimes \mathcal{O}_{E_{3}}(t)\rightarrow 0.
\end{eqnarray*}
Donc
\begin{eqnarray*}
&&\chi (E_{3},\Gamma ^{\lambda }T_{X|E_{3}}^{\ast }\otimes \mathcal{O}%
_{E_{3}}(9\left| \lambda \right| ))=\chi (E_{1},\Gamma ^{\lambda
}T_{X|E_{1}}^{\ast }\otimes \mathcal{O}_{E_{1}}(9\left| \lambda \right| )) \\
&&-\chi (E_{1},\Gamma ^{\lambda }T_{X|E_{1}}^{\ast }\otimes \mathcal{O}%
_{E_{1}}(6\left| \lambda \right| )) \\
&&=(\chi (X,\Gamma ^{\lambda }T_{X}^{\ast }\otimes \mathcal{O}_{X}(9\left|
\lambda \right| ))-\chi (X,\Gamma ^{\lambda }T_{X}^{\ast }\otimes \mathcal{O}%
_{X}(6\left| \lambda \right| ))) \\
&&-(\chi (X,\Gamma ^{\lambda }T_{X}^{\ast }\otimes \mathcal{O}_{X}(6\left|
\lambda \right| )) -\chi (X,\Gamma ^{\lambda }T_{X}^{\ast }\otimes \mathcal{O%
}_{X}(3\left| \lambda \right| ))).
\end{eqnarray*}

On termine le calcul de Riemann-Roch explicite sur Maple (cf. annexe de \cite
{Rou04}) et la proposition est d\'{e}montr\'{e}e.
\end{proof}

\bigskip

Passons maintenant \`{a} la d\'{e}monstration du th\'{e}or\`{e}me \ref{t2}:

\bigskip

\begin{proof}
Estimons maintenant
\begin{equation*}
h^{2}(X,Gr^{\bullet }E_{3,m}T_{X}^{\ast })=\underset{0\leq \gamma \leq \frac{%
m}{5}}{\sum }(\underset{\{\lambda _{1}+2\lambda _{2}+3\lambda _{3}=m-\gamma ;%
\text{ }\lambda _{i}-\lambda _{j}\geq \gamma ,\text{ }i<j\}}{\sum }%
h^{2}(X,\Gamma ^{(\lambda _{1},\lambda _{2},\lambda _{3})}T_{X}^{\ast })).
\end{equation*}

Pour $m$ suffisamment grand $\lambda $ v\'{e}rifie $\left| \lambda \right|
=\sum \lambda _{i}>4(d-5)+18.$ En effet
\begin{equation*}
\frac{4m}{5}\leq m-\gamma =\lambda _{1}+2\lambda _{2}+3\lambda _{3}\leq
6\lambda _{1}
\end{equation*}
donc
\begin{equation*}
\left| \lambda \right| \geq \lambda _{1}\geq \frac{2m}{15}.
\end{equation*}
On applique la proposition \ref{p2} et par sommation:
\begin{equation*}
h^{2}(X,Gr^{\bullet }E_{3,m}T_{X}^{\ast })\leq d(d+13)\underset{0\leq \gamma
\leq \frac{m}{5}}{\sum }(\underset{\{\lambda _{1}+2\lambda _{2}+3\lambda
_{3}=m-\gamma ;\text{ }\lambda _{i}-\lambda _{j}\geq \gamma ,\text{ }i<j\}}{%
\sum }g(\lambda ))+O(m^{8}).
\end{equation*}
Remarquons qu'\`{a} priori la sommation se fait pour $\gamma >0$ car nos
in\'{e}galit\'{e}s supposent $\lambda _{1}>\lambda _{2}>\lambda _{3},$ mais
la sommation pour $\gamma =0$ n'influence pas le terme dominant, c'est un $%
O(m^{8}).$

Il ne reste plus qu'\`{a} \'{e}valuer $\underset{0\leq \gamma \leq \frac{m}{5%
}}{\sum }(\underset{\{\lambda _{1}+2\lambda _{2}+3\lambda _{3}=m-\gamma ;%
\text{ }\lambda _{i}-\lambda _{j}\geq \gamma ,\text{ }i<j\}}{\sum }g(\lambda
))$. Ce calcul se fait par Maple:
\begin{equation*}
\underset{0\leq \gamma \leq \frac{m}{5}}{\sum }(\underset{\{\lambda
_{1}+2\lambda _{2}+3\lambda _{3}=m-\gamma ;\text{ }\lambda _{i}-\lambda
_{j}\geq \gamma ,\text{ }i<j\}}{\sum }g(\lambda ))\underset{m\rightarrow
+\infty }{\sim }\frac{49403}{252.10^{7}}m^{9}.
\end{equation*}
Et le th\'{e}or\`{e}me est d\'{e}montr\'{e}.
\end{proof}

\bigskip

On peut maintenant montrer le th\'{e}or\`{e}me \ref{teo2}:

\bigskip

\begin{proof}
On a
\begin{equation*}
h^{0}(X,E_{3,m}T_{X}^{\ast })+h^{2}(X,E_{3,m}T_{X}^{\ast })\geq \chi
(X,E_{3,m}T_{X}^{\ast })
\end{equation*}
et par \cite{Rou05}:
\begin{equation*}
\chi (X,E_{3,m}T_{X}^{\ast })=\frac{m^{9}}{81648\times 10^{6}}%
d(389d^{3}-20739d^{2}+185559d-358873)+O(m^{8}).
\end{equation*}
Par ailleurs
\begin{equation*}
h^{2}(X,E_{3,m}T_{X}^{\ast })\leq h^{2}(X,Gr^{\bullet }E_{3,m}T_{X}^{\ast
})\leq Cd(d+13)m^{9}+O(m^{8})
\end{equation*}
donc
\begin{eqnarray*}
h^{0}(X_{3},\mathcal{O}_{X_{3}}(m)) &=&h^{0}(X,E_{3,m}T_{X}^{\ast }) \\
&\geq &m^{9}(\frac{1}{81648\times 10^{6}}%
d(389d^{3}-20739d^{2}+185559d-358873) \\
&&-Cd(d+13))+O(m^{8}).
\end{eqnarray*}
Il ne reste plus qu'\`{a} evaluer pour quels degr\'{e}s
\begin{equation*}
\frac{1}{81648\times 10^{6}}d(389d^{3}-20739d^{2}+185559d-358873)-Cd(d+13)
\end{equation*}
est positif. Cela se fait par Maple. On obtient alors que $\mathcal{O}%
_{X_{3}}(m)$ est ''big'' pour $d\geq 97$ donc pour $m$ suffisamment grand:
\begin{equation*}
H^{0}(X_{3},\mathcal{O}_{X_{3}}(m)\otimes \pi _{3}^{\ast }A^{-1})\simeq
H^{0}(X,E_{3,m}T_{X}^{\ast }\otimes A^{-1})\neq 0.
\end{equation*}
\end{proof}

\section{Le cas logarithmique}

Nous allons montrer le th\'{e}or\`{e}me

\begin{theorem}
\label{t3}Soit $(\mathbb{P}^{3},X)$ vari\'{e}t\'{e} logarithmique o\`{u} X
est une hypersurface lisse de degr\'{e} $d$ de $\mathbb{P}^{3}$, alors
\begin{equation*}
h^{2}(\mathbb{P}^{3},Gr^{\bullet }E_{3,m}\overline{T_{\mathbb{P}^{3}}}^{\ast
})\leq C(d+14)m^{9}+O(m^{8})
\end{equation*}
o\`{u} C est une constante.
\end{theorem}

Montrons tout d'abord la proposition

\begin{proposition}
\label{p3}Soit $\lambda =(\lambda _{1},\lambda _{2},\lambda _{3})$ une
partition telle que $\lambda _{1}>\lambda _{2}>\lambda _{3}$ et $\left|
\lambda \right| =\sum \lambda _{i}>3d+2.$ Alors :
\begin{equation*}
h^{2}(Fl(\overline{T_{\mathbb{P}^{3}}}^{\ast }),\mathcal{L}^{\lambda
})=h^{2}(\mathbb{P}^{3},\Gamma ^{\lambda }\overline{T_{\mathbb{P}^{3}}}%
^{\ast })\leq g(\lambda )(d+14)+r(\lambda )
\end{equation*}
o\`{u} $g(\lambda )=\frac{3\left| \lambda \right| ^{3}}{2}\underset{\lambda
_{i}>\lambda _{j}}{\prod }(\lambda _{i}-\lambda _{j})$ et de plus $r$ est un
polyn\^{o}me en $\lambda $ de composantes homog\`{e}nes de plus haut
degr\'{e} 5.
\end{proposition}

\begin{proof}
On a
\begin{equation*}
\mathcal{L}^{\lambda }=(\mathcal{L}^{\lambda }\otimes \pi ^{\ast }\mathcal{O}%
_{\mathbb{P}^{3}}(3\left| \lambda \right| ))\otimes (\pi ^{\ast }\mathcal{O}%
_{\mathbb{P}^{3}}(3\left| \lambda \right| ))^{-1}=F\otimes G^{-1},
\end{equation*}
avec $F=\mathcal{L}^{\lambda }\otimes \pi ^{\ast }\mathcal{O}_{\mathbb{P}%
^{3}}(3\left| \lambda \right| )$, $G=\pi ^{\ast }\mathcal{O}_{\mathbb{P}%
^{3}}(3\left| \lambda \right| ).$ $\mathcal{L}^{\lambda }\otimes \pi ^{\ast }%
\mathcal{O}_{\mathbb{P}^{3}}(3\left| \lambda \right| )$ est positif. En
effet, on a la propri\'{e}t\'{e} g\'{e}n\'{e}rale \cite{De87} que si $E$ est
un fibr\'{e} vectoriel semi-positif i.e $E\geq 0$ alors le fibr\'{e} en
droites correspondant $\mathcal{L(}E)^{\lambda }$ est aussi semi-positif.
Ici, $E=\overline{T_{\mathbb{P}^{3}}}^{\ast }\otimes \mathcal{O}_{\mathbb{P}%
^{3}}(2)$ est semi-positif. En effet, $\overline{T_{\mathbb{P}^{3}}}^{\ast
}\otimes \mathcal{O}_{\mathbb{P}^{3}}(2)$ est engendr\'{e} par ses sections:
les sections globales de $T_{\mathbb{P}^{3}}^{\ast }\otimes \mathcal{O}_{%
\mathbb{P}^{3}}(2)$ et celles fournies par les d\'{e}riv\'{e}es
logarithmiques des \'{e}quations locales d\'{e}finissant $X.$ On a:
\begin{equation*}
\mathcal{L(}E)^{\lambda }\simeq \mathcal{L}^{\lambda }\otimes \pi ^{\ast }%
\mathcal{O}_{\mathbb{P}^{3}}(2\left| \lambda \right| )\geq 0,
\end{equation*}
donc
\begin{equation*}
\mathcal{L}^{\lambda }\otimes \pi ^{\ast }\mathcal{O}_{\mathbb{P}%
^{3}}(3\left| \lambda \right| )>0.
\end{equation*}

Soit $Y=Fl(\overline{T_{\mathbb{P}^{3}}}^{\ast }).$ Tout d'abord montrons
que $H^{i}(Y,F)=0$ pour tout $i\geq 1$ et $\lambda $ telle que $\left|
\lambda \right| =\sum \lambda _{i}>4(d-5)+18.$ Pour cela nous utilisons le
th\'{e}or\`{e}me d'annulation de Kodaira qui stipule que pour tout fibr\'{e}
en droites $A$ ample sur une vari\'{e}t\'{e} projective $Z$ complexe $%
H^{i}(Z,K_{Z}\otimes A)=0$ pour $i>0.$ En effet, regardons \`{a} quelles
conditions
\begin{equation*}
F\otimes K_{Y}^{-1}>0.
\end{equation*}
Rappelons \cite{Ma} que
\begin{equation*}
K_{Y}=\mathcal{L}^{-(5,3,1)}\otimes \pi ^{\ast }(K_{\mathbb{P}^{3}}\otimes
\det (\overline{T_{\mathbb{P}^{3}}}^{\ast })^{\otimes 3})=\mathcal{L}%
^{-(5,3,1)}\otimes \pi ^{\ast }\mathcal{O}_{\mathbb{P}^{3}}(3d-16).
\end{equation*}
Donc
\begin{equation*}
F\otimes K_{Y}^{-1}=\mathcal{L}^{\lambda +(5,3,1)}\otimes \pi ^{\ast }%
\mathcal{O}_{X}(3\left| \lambda \right| -3d+16).
\end{equation*}
Or on a
\begin{equation*}
\mathcal{L}^{\lambda +(5,3,1)}\otimes \pi ^{\ast }\mathcal{O}_{X}(2\left|
\lambda +(5,3,1)\right| )=\mathcal{L}^{\lambda +(5,3,1)}\otimes \pi ^{\ast }%
\mathcal{O}_{X}(2\left| \lambda \right| +18)\geq 0.
\end{equation*}
Par cons\'{e}quent $F\otimes K_{Y}^{-1}>0$ si
\begin{equation*}
3\left| \lambda \right| -3d+16>2\left| \lambda \right| +18
\end{equation*}
c'est-\`{a}-dire
\begin{equation*}
\left| \lambda \right| >3d+2.
\end{equation*}
Prenons un diviseur $D=\pi ^{\ast }E_{1}\in \left| G\right| ,$ lisse et
irr\'{e}ductible. On a la suite exacte:
\begin{equation*}
0\rightarrow \mathcal{O}_{Y}(F\otimes G^{-1})\rightarrow \mathcal{O}%
_{Y}(F)\rightarrow \mathcal{O}_{D}(F)\rightarrow 0.
\end{equation*}
donc la suite exacte longue en cohomologie:
\begin{equation*}
0=H^{1}(Y,\mathcal{O}_{Y}(F))\rightarrow H^{1}(D,\mathcal{O}%
_{D}(F))\rightarrow H^{2}(Y,\mathcal{O}_{Y}(F\otimes G^{-1})\rightarrow
H^{2}(Y,\mathcal{O}_{Y}(F))=0.
\end{equation*}
Donc
\begin{equation*}
h^{2}(Y,\mathcal{O}_{Y}(F\otimes G^{-1}))=h^{1}(D,\mathcal{O}_{D}(F)).
\end{equation*}
Prenons un deuxi\`{e}me diviseur $D^{\prime }=\pi ^{\ast }E_{2}\in \left|
G\right| ,$ lisse et irr\'{e}ductible, rencontrant $D$ proprement. Soit $%
Z=D\cap D^{\prime }$, $F^{\prime }=F\otimes G$ et $E_{3}=E_{1}\cap E_{2}.$
On a la suite exacte
\begin{equation*}
0\rightarrow \mathcal{O}_{D}(F^{\prime }\otimes G^{-1})\rightarrow \mathcal{O%
}_{D}(F^{\prime })\rightarrow \mathcal{O}_{Z}(F^{\prime })\rightarrow 0.
\end{equation*}
Par adjonction
\begin{equation*}
K_{D}=(K_{Y})_{|D}\otimes \mathcal{O}_{D}(D)
\end{equation*}
donc
\begin{equation*}
F_{|D}^{\prime }\otimes K_{D}^{-1}=(F\otimes K_{Y}^{-1})_{|D}>0.
\end{equation*}
Ainsi
\begin{equation*}
h^{1}(D,\mathcal{O}_{D}(F))\leq h^{0}(Z,\mathcal{O}_{Z}(F^{\prime
}))=h^{0}(Z,\mathcal{O}_{Z}(F\otimes G))\leq h^{0}(Z,\mathcal{O}%
_{Z}(F\otimes G^{2})).
\end{equation*}
Or comme pr\'{e}c\'{e}demment
\begin{equation*}
\mathcal{O}_{Z}(F\otimes G^{2})\otimes K_{Z}^{-1}=(F\otimes
K_{Y}^{-1})_{|Z}>0
\end{equation*}
donc
\begin{equation*}
h^{0}(Z,\mathcal{O}_{Z}(F\otimes G^{2}))=\chi (Z,\mathcal{O}_{Z}(F\otimes
G^{2})).
\end{equation*}
On a
\begin{equation*}
\chi (Z,\mathcal{O}_{Z}(F\otimes G^{2}))=\chi (E_{3},\Gamma ^{\lambda }%
\overline{T_{\mathbb{P}^{3}}}_{|E_{3}}^{\ast }\otimes \mathcal{O}%
_{E_{3}}(9\left| \lambda \right| )).
\end{equation*}
Par Riemann-Roch et la proposition \ref{p10}, on sait explicitement calculer
\begin{equation*}
\chi (\mathbb{P}^{3},\Gamma ^{\lambda }\overline{T_{\mathbb{P}^{3}}}^{\ast
}\otimes \otimes \mathcal{O}_{\mathbb{P}^{3}}(t)).
\end{equation*}
On a les suites exactes
\begin{eqnarray*}
(1)\text{ }0 &\rightarrow &\Gamma ^{\lambda }\overline{T_{\mathbb{P}^{3}}}%
^{\ast }\otimes \mathcal{O}_{\mathbb{P}^{3}}(t-E_{1})\rightarrow \Gamma
^{\lambda }\overline{T_{\mathbb{P}^{3}}}^{\ast }\otimes \mathcal{O}%
_{X}(t)\rightarrow \Gamma ^{\lambda }\overline{T_{\mathbb{P}^{3}}}%
_{|E_{1}}^{\ast }\otimes \mathcal{O}_{E_{1}}(t)\rightarrow 0 \\
(2)\text{ }0 &\rightarrow &\Gamma ^{\lambda }\overline{T_{\mathbb{P}^{3}}}%
_{|E_{1}}^{\ast }\otimes \mathcal{O}_{E}(t-E_{3})\rightarrow \Gamma
^{\lambda }\overline{T_{\mathbb{P}^{3}}}_{|E_{1}}^{\ast }\otimes \mathcal{O}%
_{E_{1}}(t)\rightarrow \Gamma ^{\lambda }\overline{T_{\mathbb{P}^{3}}}%
_{|E_{3}}^{\ast }\otimes \mathcal{O}_{E_{3}}(t)\rightarrow 0.
\end{eqnarray*}
Donc
\begin{eqnarray*}
&&\chi (E_{3},\Gamma ^{\lambda }\overline{T_{\mathbb{P}^{3}}}_{|E_{3}}^{\ast
}\otimes \mathcal{O}_{E_{3}}(9\left| \lambda \right| ))=\chi (E_{1},\Gamma
^{\lambda }\overline{T_{\mathbb{P}^{3}}}_{|E_{1}}^{\ast }\otimes \mathcal{O}%
_{E_{1}}(9\left| \lambda \right| )) \\
&&-\chi (E_{1},\Gamma ^{\lambda }\overline{T_{\mathbb{P}^{3}}}%
_{|E_{1}}^{\ast }\otimes \mathcal{O}_{E_{1}}(6\left| \lambda \right| )) \\
&&=(\chi (\mathbb{P}^{3},\Gamma ^{\lambda }\overline{T_{\mathbb{P}^{3}}}%
^{\ast }\otimes \mathcal{O}_{X}(9\left| \lambda \right| ))-\chi (\mathbb{P}%
^{3},\Gamma ^{\lambda }\overline{T_{\mathbb{P}^{3}}}^{\ast }\otimes \mathcal{%
O}_{\mathbb{P}^{3}}(6\left| \lambda \right| ))) \\
&&-(\chi (\mathbb{P}^{3},\Gamma ^{\lambda }\overline{T_{\mathbb{P}^{3}}}%
^{\ast }\otimes \mathcal{O}_{X}(6\left| \lambda \right| ))-\chi (\mathbb{P}%
^{3},\Gamma ^{\lambda }\overline{T_{\mathbb{P}^{3}}}^{\ast }\otimes \mathcal{%
O}_{\mathbb{P}^{3}}(3\left| \lambda \right| ))).
\end{eqnarray*}

On termine le calcul de Riemann-Roch explicite (cf. annexe de \cite{Rou04})
sur Maple et la proposition est d\'{e}montr\'{e}e.
\end{proof}

\bigskip

Passons maintenant \`{a} la d\'{e}monstration du th\'{e}or\`{e}me \ref{t3}:

\bigskip

\begin{proof}
Estimons maintenant
\begin{equation*}
h^{2}(\mathbb{P}^{3},Gr^{\bullet }E_{3,m}\overline{T_{\mathbb{P}^{3}}}^{\ast
})=\underset{0\leq \gamma \leq \frac{m}{5}}{\sum }(\underset{\{\lambda
_{1}+2\lambda _{2}+3\lambda _{3}=m-\gamma ;\text{ }\lambda _{i}-\lambda
_{j}\geq \gamma ,\text{ }i<j\}}{\sum }h^{2}(\mathbb{P}^{3},\Gamma ^{(\lambda
_{1},\lambda _{2},\lambda _{3})}\overline{T_{\mathbb{P}^{3}}}^{\ast })).
\end{equation*}

Pour $m$ suffisamment grand $\lambda $ v\'{e}rifie $\left| \lambda \right|
=\sum \lambda _{i}>3d+2.$ En effet
\begin{equation*}
\frac{4m}{5}\leq m-\gamma =\lambda _{1}+2\lambda _{2}+3\lambda _{3}\leq
6\lambda _{1}
\end{equation*}
donc
\begin{equation*}
\left| \lambda \right| \geq \lambda _{1}\geq \frac{2m}{15}.
\end{equation*}
On applique la proposition \ref{p3} et par sommation:
\begin{equation*}
h^{2}(\mathbb{P}^{3},Gr^{\bullet }E_{3,m}\overline{T_{\mathbb{P}^{3}}}^{\ast
})\leq (d+14)\underset{0\leq \gamma \leq \frac{m}{5}}{\sum }(\underset{%
\{\lambda _{1}+2\lambda _{2}+3\lambda _{3}=m-\gamma ;\text{ }\lambda
_{i}-\lambda _{j}\geq \gamma ,\text{ }i<j\}}{\sum }g(\lambda ))+O(m^{8}).
\end{equation*}
Remarquons qu'\`{a} priori la sommation se fait pour $\gamma >0$ car nos
in\'{e}galit\'{e}s supposent $\lambda _{1}>\lambda _{2}>\lambda _{3},$ mais
la sommation pour $\gamma =0$ n'influence pas le terme dominant, c'est un $%
O(m^{8}).$

Il ne reste plus qu'\`{a} \'{e}valuer $\underset{0\leq \gamma \leq \frac{m}{5%
}}{\sum }(\underset{\{\lambda _{1}+2\lambda _{2}+3\lambda _{3}=m-\gamma ;%
\text{ }\lambda _{i}-\lambda _{j}\geq \gamma ,\text{ }i<j\}}{\sum }g(\lambda
))$. Ce calcul se fait par Maple:
\begin{equation*}
\underset{0\leq \gamma \leq \frac{m}{5}}{\sum }(\underset{\{\lambda
_{1}+2\lambda _{2}+3\lambda _{3}=m-\gamma ;\text{ }\lambda _{i}-\lambda
_{j}\geq \gamma ,\text{ }i<j\}}{\sum }g(\lambda ))\underset{m\rightarrow
+\infty }{\sim }\frac{49403}{252.10^{7}}m^{9}.
\end{equation*}
Et le th\'{e}or\`{e}me est d\'{e}montr\'{e}.
\end{proof}

\bigskip

On montre alors le th\'{e}or\`{e}me \ref{teo3}:

\bigskip

\begin{proof}
On a
\begin{equation*}
h^{0}(\mathbb{P}^{3},E_{3,m}\overline{T_{\mathbb{P}^{3}}}^{\ast })+h^{2}(%
\mathbb{P}^{3},E_{3,m}\overline{T_{\mathbb{P}^{3}}}^{\ast })\geq \chi (%
\mathbb{P}^{3},E_{3,m}\overline{T_{\mathbb{P}^{3}}}^{\ast })
\end{equation*}
et par \cite{Rou05}:
\begin{eqnarray*}
\chi (\mathbb{P}^{3},E_{3,m}\overline{T_{\mathbb{P}^{3}}}^{\ast })&=&m^{9}(%
\frac{389}{81648000000}d^{3}-\frac{6913}{34020000000}d^{2} \\
&+&\frac{6299}{4252500000}d-\frac{1513}{63787500})+O(m^{8}).
\end{eqnarray*}
Par ailleurs
\begin{equation*}
h^{2}(\mathbb{P}^{3},E_{3,m}\overline{T_{\mathbb{P}^{3}}}^{\ast })\leq h^{2}(%
\mathbb{P}^{3},Gr^{\bullet }E_{3,m}\overline{T_{\mathbb{P}^{3}}}^{\ast
})\leq C(d+14)m^{9}+O(m^{8})
\end{equation*}
donc
\begin{eqnarray*}
h^{0}(\mathbb{P}^{3},E_{3,m}\overline{T_{\mathbb{P}^{3}}}^{\ast
})&\geq&
m^{9}(\frac{389}{81648000000}d^{3}-\frac{6913}{34020000000}d^{2}+\frac{6299}{%
4252500000}d \\
&-&\frac{1513}{63787500})-m^{9}C(d+14)+O(m^{8}).
\end{eqnarray*}
Il ne reste plus qu'\`{a} evaluer pour quels degr\'{e}s
\begin{equation*}
\frac{389}{81648000000}d^{3}-\frac{6913}{34020000000}d^{2}+\frac{6299}{%
4252500000}d-\frac{1513}{63787500}-C(d+14)
\end{equation*}
est positif. Cela se fait par Maple. On conclut de la m\^{e}me mani\`{e}re
que dans le cas compact:
\begin{equation*}
\text{pour }d\geq 92\text{ et }m\text{ suffisamment grand}:H^{0}(X,E_{3,m}%
\overline{T_{\mathbb{P}^{3}}}^{\ast }\otimes A^{-1})\neq 0.
\end{equation*}
\end{proof}

\end{document}